\documentclass[12pt]{article}
\usepackage{setspace}
\usepackage{amsmath,amssymb,amsthm}
\usepackage{mathrsfs}
\usepackage{geometry}
\usepackage{enumitem}
\usepackage{environ}
\usepackage{setspace}
\usepackage{titlesec}
\usepackage[all,cmtip]{xy}
\usepackage{tikz}
\usetikzlibrary{matrix,arrows}
\usepackage[colorlinks=true]{hyperref}

\NewEnviron{prf}[1][]{\begin{proof}[\bf #1Proof]\BODY\end{proof}}{}
\NewEnviron{slt}[1][]{\begin{proof}[\bf #1	解]\BODY\end{proof}}{}
\newtheorem{definition}{Definition}[section]
\newtheorem{theorem}[definition]{Theorem}
\newtheorem{lemma}[definition]{Lemma}
\newtheorem{proposition}[definition]{Proposition}
\newtheorem{remark}[definition]{Remark}
\newtheorem{example}[definition]{Example}
\newtheorem{corollary}[definition]{Corollary}

\newtheorem{question}[definition]{Question}

\newcommand\MSC[1]{\textbf{MSC(2020)}: #1}
\newcommand\address[1]{#1}
\newcommand\email[1]{\emph{Email address}: #1}
  
\setlist[enumerate,1]{label=(\roman*)}
\setlist[enumerate,2]{label=(\alph*)}

\title{\textbf{Height arguments toward the dynamical Mordell--Lang problem in arbitrary characteristic}}
\author{Junyi Xie\quad and\quad She Yang}
\date{}

\begin{document}
\begin{spacing}{1.25}

\maketitle

\begin{abstract}
We use height arguments to prove two results about the dynamical Mordell--Lang problem.

(i) For an endomorphism of a projective variety, the return set of a dense orbit into a curve is finite if any cohomological Lyapunov multiplier of any iteration is not an integer.

(ii) Let $f\times g:X\times C\rightarrow X\times C$ be an endomorphism, where $f$ and $g$ are surjective endomorphisms of a projective variety $X$ and a projective curve $C$, respectively. If the degree of $g$ is greater than the first dynamical degree of $f$, then the return sets of the system $(X\times C,f\times g)$ have the same form as the return sets of the system $(X,f)$.

Using the second result, we deal with the case of split self-maps of products of curves, for which the degrees of the factors are pairwise distinct.

In the cases that the height argument cannot be applied, we find examples that show the return set can be very complicated --- more complicated than experts once imagined --- even for endomorphisms of tori with zero entropy. One may compare them with the conjectures and known results stated in the works of Corvaja--Ghioca--Scanlon--Zannier [6] and the authors [36].
\end{abstract}

\MSC{Primary: 37P55, 37P30; Secondary: 37P15.}

\section{Introduction}

In this paper, we work over an algebraically closed field $K$ of arbitrary characteristic. Our original purpose was to prove some results in positive characteristic, but our theorems turn out to also be valid in zero characteristic.

Unless otherwise specified, the varieties and maps are over $K$. As a matter of convention, every variety is assumed to be integral, but the closed subvarieties can be reducible. For a rational map $f:X\dashrightarrow Y$ between two varieties, we denote $\mathrm{Dom}(f)\subseteq X$ as the domain of definition of $f$. Let $X$ be a variety and let $f$ be a rational self-map of $X$. For a point $x\in X(K)$, we say the orbit $\mathcal{O}_{f}(x):=\{f^{n}(x)|\ n\in\mathbb{N}\}$ is well-defined if every iterate $f^{n}(x)$ lies in $\text{Dom}(f)$. We denote $\mathbb{N}=\mathbb{Z}_{+}\cup\{0\}$. An arithmetic progression is a set of the form $\{mk+l|\ k\in\mathbb{Z}\}$ for some $m,l\in\mathbb{Z}$ and an arithmetic progression in $\mathbb{N}$ is a set of the form $\{mk+l|\ k\in\mathbb{N}\}$ for some $m,l\in\mathbb{N}$.

The dynamical Mordell--Lang conjecture is one of the core problems in the field of arithmetic dynamics. It asserts that for any rational self-map $f$ of a variety $X$ over $\mathbb{C}$, the return set $\{n\in\mathbb{N}|\ f^{n}(x)\in V(\mathbb{C})\}$ is a finite union of arithmetic progressions in $\mathbb{N}$, where $x\in X(\mathbb{C})$ is a point such that the orbit $\mathcal{O}_{f}(x)$ is well-defined and $V\subseteq X$ is a closed subvariety. There is an extensive literature on various cases of this 0-DML conjecture (``0" stands for the characteristic of the base field). Two significant cases are as follows.
\begin{enumerate}
\item
If $X$ is a quasi-projective variety over $\mathbb{C}$ and $f$ is an \'etale endomorphism of $X$, then the 0-DML conjecture holds for $(X,f)$. See \cite{Bel06} and \cite[Theorem 1.3]{BGT10}.
\item
If $X=\mathbb{A}_{\mathbb{C}}^{2}$ and $f$ is an endomorphism of $X$, then the 0-DML conjecture holds for $(X,f)$. See \cite{Xie17} and \cite[Theorem 4]{Xie}.
\end{enumerate}

One can consult \cite{BGT16,Xie} and the references therein for further known results.

The statement of the 0-DML conjecture fails when the base field has positive characteristic. See \cite[Example 3.4.5.1]{BGT16} for an example. Indeed, the return set can be very complicated in positive characteristic. Ghioca and Scanlon once proposed a $p$DML conjecture on the form of the return set. However, as showed in \cite{XY}, the return set can be even more complicated than what they conjectured. In the last section of this paper, we give examples of return sets having an unprecedented form. Such a form is even beyond the scope of ``widely $p$-normal sets" as defined in \cite[Definition 1.1]{XY}. See Question 5.9 for some discussion.

The $p$DML problem is known to be \emph{very} hard. Indeed, it is proved in \cite{CGSZ21} that the $p$DML problem for endomorphisms of tori is \emph{equivalent} to solving certain hard Diophantine equations. Consequently, not much is known about this $p$DML problem. One can consult \cite{CGSZ21}, \cite[Theorem 1.4, Theorem 1.5]{Xie23}, \cite{Yang24}, and \cite{XY} for references.

~

In this paper, we use height arguments to study certain cases of the DML problem in arbitrary characteristic.

Before stating our main theorems, we first recall the definitions of the dynamical degrees and the cohomological Lyapunov multipliers of an algebraic dynamical system. We only state the definitions for endomorphisms of projective varieties for simplicity, since this is the only case that we encounter in this paper. However, we note that these concepts also make sense for dominant rational self-maps.

Let $f:X\rightarrow X$ be a surjective endomorphism of a projective variety. Let $L$ be a big and nef line bundle on $X$. Then for every $i\in\{0,\dots,\text{dim}(X)\}$, the $i$-th \emph{dynamical degree} of $f$ is $\lambda_i(f)=\lim\limits_{n\rightarrow\infty}((f^n)^*L^i\cdot L^{n-i})^{\frac{1}{n}}$. These are very important quantities that measure the complexity of the algebraic dynamical system and have been carefully studied in the literature. See for example \cite{DS05,Dang20,Tru20}, and \cite[Section 2.1]{Xie23}. In particular, the limits exist and do not depend on the choice of $L$. For example, we have $\lambda_0(f)=1$ and $\lambda_{\text{dim}(X)}(f)=\text{deg}(f)$.

In the setting above, we define the $i$-th \emph{cohomological Lyapunov multiplier} $\mu_i(f)=\frac{\lambda_i(f)}{\lambda_{i-1}(f)}$ for $i\in\{1,\dots,\text{dim}(X)\}$. This concept was introduced by the first author in \cite{Xiec}, and he shows in \cite{Xieb} that it has many interesting properties. For example, let us denote $f^*:\text{N}^1(X)_{\mathbb{R}}\rightarrow\text{N}^1(X)_{\mathbb{R}}$ as the $\mathbb{R}$-linear self-map of the numerical class group induced by $f$. It is not hard to see that $\lambda_1(f)=\mu_1(f)$ is the spectral radius of $f^*$, but in fact \emph{all} of the $\mu_i(f)$ are eigenvalues of $f^*$. See \cite[Theorem 1.4]{Xieb} for this, and please consult subsection 2.1 for more information about the cohomological Lyapunov multipliers.

Now we can state our first main theorem. In this paper, we denote $\mathfrak{Root}=\{a^{\frac{1}{n}}|\ a,n\in\mathbb{Z}_+\}$. Here, by $a^{\frac{1}{n}}$, we mean the unique positive real root.

\begin{theorem}
Let $X$ be a projective variety and let $f$ be a surjective endomorphism of $X$. Suppose that $\mu_i(f)\notin\mathfrak{Root}$ for every $i\in\{1,\dots,\mathrm{dim}(X)\}$. Let $C\subseteq X$ be a closed subcurve and let $x\in X(K)$ be a point. If $\overline{\mathcal{O}_f(x)}=X$, then $\mathcal{O}_{f}(x)\cap C(K)$ is a finite set.
\end{theorem}

\begin{remark}
By tracking through the proof, one can see that the same statement holds for every irreducible closed subvariety $C\subseteq X$ of Picard number 1 if $X$ is normal. This has some interest, since the Picard number may increase when taking a subvariety. But for simplicity, we will focus on the case of curves.
\end{remark}

We make some comments on our assumption that $\mu_i(f)\notin\mathfrak{Root}$. Firstly, we remind that this condition is not expected to be truly necessary in characteristic 0, as the dynamical Mordell--Lang conjecture suggested. But in the case of characteristic $p>0$, the examples in Section 5 imply that at least we may need to assume $\log_{p}\mu_i(f)\notin\mathbb{Q}_{\geq0}$. In an ongoing work, we will discuss this problem in more detail.

By definition, an endomorphism $f$ of $X$ is \emph{cohomologically hyperbolic} if $\mu_i(f)\neq1$ for every $i\in\{1,\dots,\text{dim}(X)\}$. Since $1\in\mathfrak{Root}$, we see that our requirement above forces $f$ to be cohomologically hyperbolic. Philosophically speaking, this requirement depicts cohomologically hyperbolic endomorphisms that are opposite to the polarized ones.

We now briefly summarize the proof strategy of Theorem 1.1 to give the reader an overview and indicate how the condition $\mu_i(f)\notin\mathfrak{Root}$ is used. Firstly, by using an Albanese argument, we reduce to proving the two extremal cases, i.e. the cases that $\mathrm{Pic}^0(X)$ is trivial or $X$ is itself an abelian variety. Then we use the height argument. We explain the first case in the following. This case has the advantage that $\mathrm{Pic}(X)_{\mathbb{C}}=\mathrm{N}^1(X)_{\mathbb{C}}$.

The key point in running the height argument is finding \emph{two different} speeds of height growth. This is illustrated in Proposition 3.1. The condition (ii) given in there guarantees that we can get two different speeds of height growth of the orbit $\mathcal{O}_f(x)$ when calculated by $L_1$ and $L_2$, respectively. But condition (i) says that they can control each other on the subsequence $\mathcal{O}_f(x)\cap C(K)$. Hence the conclusion follows.

But how to find $L_1,L_2$ and $\mu_1,\mu_2$ satisfying the conditions of Proposition 3.1? This is where the cohomological Lyapunov multipliers come into play. By Theorem 2.3, we can first find one pair of the desired $(L,\mu)$, in which $\mu$ is a certain cohomological Lyapunov multiplier of $f$. Then, by using the trick of considering Galois conjugates, we can get the second pair of $(L,\mu)$ from the first pair. See Proposition 3.3. Note that we have assumed $\mathrm{Pic}(X)_{\mathbb{C}}=\mathrm{N}^1(X)_{\mathbb{C}}$ in this discussion. The assumption that $\mu_i(f)\notin\mathfrak{Root}$ is used here in order to guarantee that we can indeed get a new pair different from the original one. See Lemma 3.4.

~

Theorem 1.1 can deduce the following corollary. Please see Definition 2.4 for the definition of amplified endomorphisms. One can see that part (ii) of the corollary below is a generalization of \cite[Theorem 1.4]{Xie23}.

\begin{corollary}
Let $X$ be a projective variety and let $f$ be a surjective endomorphism of $X$.
\begin{enumerate}
\item
If $\mathrm{dim}(X)=2$ and $\lambda_1(f)\notin\mathfrak{Root}$, then $(X,f)$ satisfies the $\text{DML}_0$ property (see Definition 1.6).
\item
Suppose that $\{\mu_1(f),\dots,\mu_{\mathrm{dim}(X)}(f)\}$ has an empty intersection with the interval $[1,\mathrm{deg}(f)]$. Then the conclusion of Theorem 1.1 holds for $f$. In particular, the conclusion holds for $f$ if $f$ is a cohomologically hyperbolic automorphism or $f$ satisfies $\lambda_1(f)>\lambda_2(f)$.
\item
If $f$ is an amplified automorphism, then $\mathcal{O}_{f}(x)\cap C(K)$ is a finite set for every closed subcurve $C\subseteq X$ and every point $x\in X(K)$. We do not need to assume $\overline{\mathcal{O}_f(x)}=X$ here.
\end{enumerate}
\end{corollary}

In accordance with the referees' requests, we now try to present some explicit examples of Theorem 1.1 and Corollary 1.3. From the \emph{numerical} point of view, it seems that the assumption that $\mu_i(f)\notin\mathfrak{Root}$ of Theorem 1.1 should hold for ``general" endomorphisms. But the explicit examples of endomorphisms of projective varieties are mainly (products of) polarized systems, which are indeed philosophically opposite to the systems satisfying our condition. In particular, in characteristic 0, we find it hard to give concrete examples of Theorem 1.1 that are not covered by previous results. Therefore, we focus on the positive characteristic case and give two examples below. They do not appear to be covered by known results (in positive characteristic).
\begin{enumerate}
\item
Endomorphisms of abelian varieties. Let $G$ be an abelian variety and let $X=G^N$. Let $f$ be the endomorphism of $X$ induced by a matrix $A\in M_n(\mathbb{Z})$. By using \cite[Theorem 1.9]{Hu24} or results in \cite{Xies}, one can verify that every multiplier $\mu_i(f)$ is the square of the modulus of an eigenvalue of $A$. This gives many examples of Theorem 1.1.

\item
Cohomologically hyperbolic automorphisms of hyperk\"ahler varieties. As we mentioned above, Corollary 1.3(ii) is a generalization of \cite[Theorem 1.4]{Xie23}, which dealt with the case of surfaces. We refer to \cite[Section 4]{Ogu09} for some examples of cohomologically hyperbolic automorphisms of hyperk\"ahler varieties. In particular, according to Example 2 in loc. cit., there are some such automorphisms that are not induced from automorphisms of surfaces. Hence we get some examples in this type that are not covered by previous results.
\end{enumerate}

~

Before stating our second main result, we give some definitions regarding the possible forms of return sets. All of these seemingly strange forms are devoted for the dynamical Mordell--Lang problem in positive characteristic. The definition of $p$-normal sets was firstly introduced in \cite{Der07} for the Skolem--Mahler--Lech problem in positive characteristic, and the definition of widely $p$-normal sets was introduced in \cite{XY} for describing the return sets of bounded-degree rational self-maps.

\begin{definition}
Suppose $\mathrm{char}(K)=p>0$.

Let $q=p^{e}$ for some positive integer $e$. Suppose that $d\in\mathbb{Z}_{+},r\in\mathbb{N}$ and $c_{0},c_{ij}\in\mathbb{Q}$ where $(i,j)\in\{1,\dots,d\}\times\{0,\dots,r\}$. Then we define
$$
S_{q,d,r}(c_{0};c_{ij})=\{c_{0}+\sum\limits_{i=1}^{d}\sum\limits_{j=0}^{r} c_{ij}q^{2^{j}n_{i}}|\ n_{1},\dots,n_{d}\in\mathbb{N}\}.
$$
\begin{enumerate}
\item
We define a \emph{widely} $p$\emph{-normal set in} $\mathbb{Z}$ as a union of finitely many arithmetic progressions along with finitely many \emph{subsets of} $\mathbb{Z}$ of the form $S_{q,d,r}(c_{0};c_{ij})$ as above. A \emph{widely} $p$\emph{-normal set in} $\mathbb{N}$ is a subset of $\mathbb{N}$ which is, up to a finite set, equal to the intersection of a widely $p$-normal set in $\mathbb{Z}$ and $\mathbb{N}$. To simplify the notation, we call widely $p$-normal sets in $\mathbb{N}$ as ``type 2 sets".
\item
We define a $p$\emph{-normal set} in $\mathbb{Z}$ as a union of finitely many arithmetic progressions along with finitely many sets of the form $S_{q,d,0}(\frac{c_{0}}{q-1};\frac{c_{i}}{q-1})$ as above in which $q$ is a power of $p$ and $c_0,c_1,\dots,c_d$ are integers satisfying $q-1\mid c_0+c_1+\cdots+c_d$. A $p$\emph{-normal set in} $\mathbb{N}$ is a subset of $\mathbb{N}$ which is, up to a finite set, equal to the intersection of a $p$-normal set in $\mathbb{Z}$ and $\mathbb{N}$. We call $p$-normal sets in $\mathbb{N}$ as ``type 1 sets".
\end{enumerate}
\end{definition}

Here we say two sets $S$ and $T$ are equal up to a finite set if the symmetric difference $(S\backslash T)\cup(T\backslash S)$ is finite, as in $\cite{Der07}$. Parallel with the type 1 and type 2 sets described above, a set which is a finite union of arithmetic progressions in $\mathbb{N}$ will be called as a type 0 set.

In order to eliminate possible confusions, we remark that the notions of type 1 and type 2 sets \emph{depend on} the prime $p=\text{char}(K)$. In particular, they do not make sense if $\text{char}(K)=0$. The letter ``$p$" in the word ``$p$-normal" should be regarded as an \emph{abbreviation} of $\text{char}(K)$ instead of a formal symbol.

\begin{remark}
For each $\epsilon\in\{0,1,2\}$, one can verify that the union and \emph{intersection} of two sets of type $\epsilon$ described above is also a set of type $\epsilon$.
\end{remark}

\begin{definition}
Let $X$ be a variety and let $f:X\dashrightarrow X$ be a dominant rational self-map. For each $\epsilon\in\{0,1,2\}$, we say that the dynamical system $(X,f)$ satisfies the $\text{DML}_{\epsilon}$ property if for every $x\in X(K)$ with a well-defined orbit and every closed subvariety $V\subseteq X$, the return set $\{n\in\mathbb{N}|\ f^{n}(x)\in V(K)\}$ is a set of type $\epsilon$.
\end{definition}

Notice that the definitions of $\text{DML}_1$ and $\text{DML}_2$ properties only make sense in positive characteristic. So whenever we mention them, we tacitly assume that the circumstance is of positive characteristic.

\begin{remark}
Let $\epsilon\in\{0,1,2\}$. If the dynamical system $(X,f^{n_0})$ satisfies the $\text{DML}_{\epsilon}$ property for some positive integer $n_0$, then so does $(X,f)$. However, it seems not easy to prove that the $\text{DML}_{\epsilon}$ property is stable under change of birational models (although we tend to believe that this is true).
\end{remark}

Now we can state our second main result.

\begin{theorem}
Let $X$ be a projective variety and let $C$ be a projective curve. Let $f:X\rightarrow X$ and $g:C\rightarrow C$ be surjective endomorphisms. Suppose that $\lambda_1(f)<\mathrm{deg}(g)$ and the dynamical system $(X,f)$ satisfies the $\text{DML}_{\epsilon}$ property in which $\epsilon\in\{0,1,2\}$. Then $(X\times C,f\times g)$ also satisfies $\text{DML}_{\epsilon}$ property.
\end{theorem}

Indeed, Theorem 1.8 also holds when $f$ and $g$ are dominant rational self-maps. But for simplicity and not letting details obscure main ideas, we focus on the case of endomorphisms.

Theorem 1.8 illustrates the connection between Kawaguchi--Silverman type results and the dynamical Mordell--Lang problem for split endomorphisms of products, which we briefly explain here. The reader may refer to Proposition 4.1 and Remark 4.4 for more on this.

As we have mentioned before, we need \emph{two different} speeds of height growth. In the split setting, the two different speeds come from the two factors. In Theorem 1.8, the growth on $X$ is slower due to the Kawaguchi--Silverman--Matsuzawa upper bound, which roughly says that exponentially the height growth cannot be quicker than the exponential function with base $\lambda_1(f)$. Please see the proof of Proposition 4.1 for more information. On the contrary, the growth on $C$ should be faster, as the dynamical system on $C$ is polarized with degree $\mathrm{deg}(g)>\lambda_1(f)$. This difference of height growth forces a preperiodicity on the faster factor. See Proposition 4.1. We need the faster factor to have dimension 1 because this restricts the form of irreducible closed subvarieties $V\subseteq X\times C$: either the projection of $V$ to $X$ is generically finite onto its image, or $V$ has the form $W\times C$ for some $W\subseteq X$.

In \cite[Corollary 1.5]{SXY}, we further develop this idea and obtain another result using this philosophy.

~

We can use Theorem 1.8 to study the split rational self-maps of product of curves. We mention that in \cite[Section 5.8]{BGT16} and \cite{CGLN}, a similar height argument is used to prove certain relative unlikely intersection results. Namely, for endomorphisms $f$ and $g$ of the projective line, they concern about the relationship between $f$ and $g$ provided that there exists an infinite intersection of orbits $\mathcal{O}_f(x)\cap\mathcal{O}_g(y)$ for some starting points $x$ and $y$.

\begin{corollary}
Let $n\geq m$ be positive integers. Let $C_1,\dots,C_n$ be projective curves and let $g_1:C_1\dashrightarrow C_1,\dots,g_n:C_n\dashrightarrow C_n$ be dominant rational self-maps. Suppose $1=\mathrm{deg}(g_1)=\cdots=\mathrm{deg}(g_m)<\mathrm{deg}(g_{m+1})<\cdots<\mathrm{deg}(g_n)$. Then the dynamical system $(C_1\times\cdots\times C_n,g_1\times\cdots\times g_n)$ satisfies the $\text{DML}_2$ property. Moreover, it satisfies the $\text{DML}_0$ property if $m=1$ or if $\mathrm{char}(K)=0$.
\end{corollary}

The statement seems to require that there has to exist at least one component of degree $1$. However, this is not really needed as otherwise we may just add an identity component.

Under this situation, we can control the complexity of the return set.

\begin{remark}
Suppose $\mathrm{char}(K)=p>0$. Denote $\mathrm{pr}$ as the projection $C_1\times\cdots\times C_n\rightarrow C_1\times\cdots\times C_m$. Let $V\subseteq C_1\times\cdots\times C_n$ be a closed subvariety. By tracking through the proof of Corollary 1.9 and taking \cite[Remark 4.9]{XY} into account, we can see that one may let all of the ``widely $p$-sets" be involved in the return set into $V$ have the form $\{c_{0}+\sum\limits_{i=1}^{d}\sum\limits_{j=0}^{r_i} c_{ij}q^{2^{j}n_{i}}|\ n_{1},\dots,n_{d}\in\mathbb{N}\}$ in which $d+\sum\limits_{i=1}^{d}r_i+|\{i|\ r_i>0,1\leq i\leq d\}|\leq\mathrm{dim}(\mathrm{pr}(V))$.
\end{remark}

~

Now we compare our methods and results with what was previously known about the dynamical Mordell--Lang conjecture. In characteristic 0, the two most remarkable known cases mentioned at the beginning of this paper mainly do not use the height argument. Thus, in the following, we follow our original motivation and focus on the positive characteristic case. The two main known cases of the $p$DML problem are as follows.
\begin{enumerate}
\item
In \cite{CGSZ21} and \cite{XY}, some endomorphisms of semiabelian varieties are studied. In particular, the $p$DML problem for translations of algebraic groups is settled in \cite{XY}. These works rely on the Mordell--Lang type results in positive characteristic proved by Hrushovski \cite{Hru96} and Moosa--Scanlon \cite{MS04}.

\item
In \cite{Yang24}, the second author dealt with certain endomorphisms of the projective space. These morphisms are called ``totally inseparable liftings of Frobenius" and indeed are specific things in positive characteristic. There are no zero characteristic analogues of them. The proof uses jet spaces.
\end{enumerate}

In this paper, we try to apply height arguments toward the dynamical Mordell--Lang conjecture in arbitrary characteristic. The idea of using heights to solve problems in arithmetic dynamics dates back to almost 20 years ago. For example, see \cite{BGKT10}, \cite[Theorem 8.1]{GTZ12}, and \cite[Lemma 7.24]{MS14}. For previous approaches using height arguments to study the $p$DML problem, see \cite{Xie14}, \cite[Theorem 1.4]{Xie23}, and \cite{Nel}. Our Theorem 1.1 and Corollary 1.3 can be regarded as generalizations of \cite[Theorem 1.4]{Xie23}, which dealt with cohomologically hyperbolic automorphisms of projective surfaces. Our proof needs some more advanced knowledge about the cohomological Lyapunov multipliers discovered by the first author. See subsection 2.1. In comparison, Theorem 1.8 and Corollary 1.9 are more elementary and explicit as the ``two different speeds of growth" automatically appear in such cases. The spirit of them is related to the work in \cite{Nel}, which dealt with certain endomorphisms of tori in positive characteristic. 

Notice that the results of \cite[Theorem 1.4]{Xie23} and \cite{Nel} are known in characteristic 0 before. On the contrary, our results mentioned above seem to be new in characteristic 0. We make a remark now, which claims that we shall focus on the positive characteristic.

\begin{remark}(About the characteristic)
Starting from Section 3, we will assume that our base field $K$ has a positive characteristic. We need a height notion on function fields that satisfies the Northcott property. In the positive characteristic case, the Weil height machinery is adequate for our purpose. But in zero characteristic, we need to use the more delicate notion of Moriwaki's height \cite{Mor00}. We shall briefly review the definition and basic properties of this notion in subsection 2.3. One can see that this machinery has all good properties that the Weil height machinery has, and the careful reader can verify that our argument is still valid line by line in the zero characteristic case by using Moriwaki's height.
\end{remark}

~

This paper also contains a section with examples for which the height argument cannot be applied, and hence the return sets can be \emph{very} complicated (in positive characteristic). Here, we only state a result which will be proved in subsection 5.2. This statement disproves \cite[Conjecture 10.2]{Der07}, which asserts that automorphisms of affine varieties satisfy the $\text{DML}_1$ property.

\begin{proposition}
There is an automorphism of a linear torus with zero entropy which does not satisfy the $\text{DML}_1$ property.
\end{proposition}

\begin{remark}
\begin{enumerate}
\item
The heuristic Example 5.6 suggests that there should exist endomorphisms of linear tori with zero entropy which do not satisfy the $\text{DML}_2$ property. So there are new phenomenon taking place in this type of endomorphisms, beyond the scope of \cite{XY}.

\item
Assume Vojta's conjecture. Then \cite[Theorem 1.6]{CGSZ21} says that the endomorphisms of tori satisfy the $\text{DML}_1$ property. However, we do \emph{not} think that we have disproved Vojta's conjecture. See subsection 5.2 for more information.
\end{enumerate}
\end{remark}

At the end of the Introduction, we outline the structure of this paper. In Section 2, we make some preparations about the cohomological Lyapunov multipliers and the Weil height machinery. Then in Section 3, we prove Theorem 1.1 and the corollaries about cohomologically hyperbolic endomorphisms. We prove Theorem 1.8 and Corollary 1.9 in Section 4. Finally, in Section 5, we propose various examples in positive characteristic for which the return sets are complicated. There are two origins of the examples --- the Frobenius and the zero entropy endomorphisms of algebraic groups.

\section{Preparations}

In this section, we make some technical preparations. In subsection 2.1, we recall some knowledge about the cohomological Lyapunov multipliers obtained in \cite{Xieb}. Then we review the Weil height machinery in subsection 2.2. In subsection 2.3, we recall the Moriwaki height machinery for arithmetic function fields. We shall use Weil's height machinery for finitely generated fields over $\mathbb{F}_p$. For finitely generated fields over $\mathbb{Q}$, Moriwaki's height machinery is a perfect replacement of the Weil height machinery.

\subsection{Cohomological Lyapunov multipliers}

In this subsection, we fix a projective variety $X$ and a surjective endomorphism $f$ of $X$. We have defined the cohomological Lyapunov multipliers $\mu_1(f),\dots,\mu_{\text{dim}(X)}(f)$ in the Introduction. Now we recall some of their properties.

We start with the proposition below. The first four parts of the proposition are immediate consequences of basic properties of the dynamical degrees, while the fifth part is a corollary of the relative dynamical degree formula. See \cite[Theorem 1.3(2)]{Tru20}.

\begin{proposition}
Let $Y$ be a projective variety and let $g$ be a surjective endomorphism of $Y$. In parts (iv)(v), let $\pi:X\rightarrow Y$ be a surjective morphism such that $\pi\circ f=g\circ\pi$. The sets in the following statements should be comprehended as multiple-sets.
\begin{enumerate}
\item
We have $\mu_1(f)\geq\cdots\geq\mu_{\mathrm{dim}(X)}(f)>0$.
\item
We have $\mu_i(f^n)=\mu_i(f)^n$ for every $i\in\{1,\dots,\mathrm{dim}(X)\}$ and every positive integer $n$.
\item
Let $f\times g$ be the surjective endomorphism of $X\times Y$ induced by $f$ and $g$. Then we have $\{\mu_1(f\times g),\dots,\mu_{\mathrm{dim}(X)+\mathrm{dim}(Y)}(f\times g)\}=\{\mu_1(f),\dots,\mu_{\mathrm{dim}(X)}(f)\}\sqcup\{\mu_1(g),\dots,\mu_{\mathrm{dim}(Y)}(g)\}$.
\item
If $\pi$ is generically finite, then we have $\{\mu_1(f),\dots,\mu_{\mathrm{dim}(X)}(f)\}=\{\mu_1(g),\dots,\mu_{\mathrm{dim}(Y)}(g)\}$.
\item
If $Y$ is smooth, then we have $\{\mu_1(g),\dots,\mu_{\mathrm{dim}(Y)}(g)\}\subseteq\{\mu_1(f),\dots,\mu_{\mathrm{dim}(X)}(f)\}$.
\end{enumerate}
\end{proposition}

Next, we recall the main properties of the cohomological Lyapunov multipliers that will be used in this article. For a projective variety $X$, we denote $\mathrm{N}^1(X)$ as the numerical class group of line bundles on $X$. The theorem of the base guarantees that $\mathrm{N}^1(X)$ is a finite free $\mathbb{Z}$-module. We denote $f^*:\text{N}^1(X)_{\mathbb{R}}\rightarrow\text{N}^1(X)_{\mathbb{R}}$ as the pull-back map induced by $f$.

\begin{theorem}(\cite[Theorem 1.4]{Xieb})
We have
$$
\{\mu_1(f),\dots,\mu_{\mathrm{dim}(X)}(f)\}=\{\alpha\in\mathbb{R}|\ \mathrm{Im}(f^*-\alpha)\cap\mathrm{Big}(X)=\emptyset\},
$$
in which $\mathrm{Big}(X)\subseteq\mathrm{N}^1(X)_{\mathbb{R}}$ is the big cone of $X$. In particular, all of the $\mu_i(f)$ are eigenvalues of $f^*$. Hence they are algebraic integers.
\end{theorem}

The theorem below is a combination of Theorem 2.2 and \cite[Theorem 1.3]{Xieb}.

\begin{theorem}
Let the notations be as above. Then the linear subspace
$$
\sum\limits_{i=1}^{\mathrm{dim}(X)}\mathrm{ker}(f^*-\mu_i(f))^{\rho}\subseteq\mathrm{N}^1(X)_{\mathbb{R}}
$$
has a nonempty intersection with $\mathrm{Big}(X)$, in which $\rho$ is the rank of $\mathrm{N}^1(X)$.
\end{theorem}

Now we introduce the definitions of \emph{cohomologically hyperbolic} and \emph{amplified} endomorphisms.

\begin{definition}
\begin{enumerate}
\item
The endomorphism $f$ is said to be \emph{cohomologically hyperbolic} if
$$
1\notin\{\mu_1(f),\dots,\mu_{\mathrm{dim}(X)}(f)\}.
$$
\item
The endomorphism $f$ (of the \emph{projective} variety $X$) is said to be \emph{amplified} if there exists a line bundle $L$ on $X$ such that $f^*L-L$ is ample.
\end{enumerate}
\end{definition}

The connection between these two types of morphisms is revealed by the following theorem. Here, by a ``subsystem" of an iteration $f^n$, we mean an $f^n$-invariant irreducible closed subvariety of $X$ endowed with the induced surjective endomorphism.

\begin{theorem}(\cite[Theorem 1.5]{Xieb})
The endomorphism $f$ is amplified if and only if every subsystem of every iteration of $f$ is cohomologically hyperbolic. In particular, amplified endomorphisms are cohomologically hyperbolic.
\end{theorem}

We also need the following result which says that the set of periodic points of a cohomologically hyperbolic endomorphism is dense.

\begin{theorem}(\cite[Theorem 1.12]{Xiec})
If $f$ is cohomologically hyperbolic, then the set of $f$-periodic closed points is dense in $X$.
\end{theorem}

\subsection{The Weil height machinery}

In this subsection, we recall the Weil height machinery following \cite[Chapter 2]{Ser97}. In the following, we will always let $k$ be a field of positive characteristic.

\begin{definition}
Let $M_k$ be a family of non-archimedean discrete absolute values on $k$. Then every $v\in M_k$ has the form $|x|_v=c^{-v(x)}$ where $v:k\rightarrow\mathbb{Z}\cup\{\infty\}$ is a discrete valuation and $c>1$. Suppose
\begin{enumerate}
\item
for all $x\in k^{\times}$ one has $|x|_v=1$ for all but finitely many $v\in M_k$, and
\item
for all $x\in k^{\times}$ one has $\prod\limits_{v\in M_k}|x|_v=1$.
\end{enumerate}
Then we say that $k$ is \emph{equipped with a product formula}.
\end{definition}

If $k$ is a field with a product formula, then we can define the naive logarithmic height function $h$ on the projective space $\mathbb{P}^{N}(k)$ in the usual way.

\begin{definition}
We say a product formula field has the \emph{Northcott property} if $\{x\in k|\ h(x)\leq A\}$ is a finite set for every $A>0$.
\end{definition}

The following statement should be well-known to experts, but we will sketch a proof because of lack of reference.

\begin{proposition}
Let $k$ be a finitely generated field extension of $\mathbb{F}_p$ of positive transcendence degree. Then we can make $k$ into a product formula field which satisfies the Northcott property.
\end{proposition}

\begin{prf}
Let $\{t_1,\dots,t_n\}$ be a separable transcendence basis of $k/\mathbb{F}_p$. Denote $k_0=\mathbb{F}_p(t_1,\dots,t_n)$. Then $k/k_0$ is a finite separable extension. By extending $k$, we may assume that $k$ is a finite Galois extension of $k_0$ without loss of generality. As we have fixed a set of transcendence basis of $k_0/\mathbb{F}_p$, there is a natural way to give a product formula on $k_0$. Then $k_0$ is a product formula field with the Northcott property because its constant field is finite. Let $M_{k_0}$ be the family of absolute values on $k_0$. We shall construct the family of absolute values $M_k$ by extending the absolute values in $M_{k_0}$.

We recall two basic facts about the extension of absolute values to finite Galois extensions. They can be easily deduced from the knowledge of $\cite[\text{Chapter}\ 1]{Lan83}$. Let $v\in M_{k_0}$ be a non-archimedean discrete absolute value. Then the following statements hold.
\begin{enumerate}
\item
There are only finitely many extension of absolute values of $v$ on $k$. Let them be $w_1,\dots,w_{g(v)}$, in which $g(v)$ is the finite number of extensions. Then each $w_i$ is a non-archimedean discrete absolute value on $k$.
\item
The action of $\mathrm{Gal}(k/k_0)$ on $\{w_1,\dots,w_{g(v)}\}$ given by $|x|_{\sigma(w)}=|\sigma^{-1}(x)|_w$ is transitive.
\end{enumerate}

For each $w_i$, we define a non-archimedean discrete absolute value $||\cdot||_{w_i}$ on $k$ by $||x||_{w_i}=|x|_{w_i}^{\frac{1}{g(v)}}$. Let $M_k$ be the family $\{||\cdot||_w:w\mid v\text{ for some }v\in M_{k_0}\}$. Then one can see that the family $M_k$ satisfies condition (i) in Definition 2.7. Moreover, we calculate that $\prod\limits_{w\mid v}||x||_w=\prod\limits_{w\mid v}|x|_{w}^{\frac{1}{g(v)}}=\prod\limits_{\sigma\in\mathrm{Gal}(k/k_0)}|x|_{\sigma(w_1)}^{\frac{1}{[k:k_0]}}=\prod\limits_{\sigma\in\mathrm{Gal}(k/k_0)}|\sigma(x)|_{w_1}^{\frac{1}{[k:k_0]}}=|N_{k/k_0}(x)|_{w_1}^{\frac{1}{[k:k_0]}}=|N_{k/k_0}(x)|_{v}^{\frac{1}{[k:k_0]}}$ for every $v\in M_{k_0}$ and every $x\in k$ and  then see (ii) also holds.

We have seen that the family $M_k$ equip a product formula on $k$. Now we show that the product formula field $k$ satisfies the Northcott property. We need the following facts.
\begin{enumerate}
\item
Height is invariant under Galois conjugate. This is because the multiple-set $\{||\sigma(x)||_w:v\mid w\}=\{||x||_{\sigma^{-1}(w)}:v\mid w\}=\{||x||_{w}:v\mid w\}$ for every $\sigma\in\mathrm{Gal}(k/k_0)$, every $x\in k$ and every $v\in M_{k_0}$.
\item
If $x\in k_0$, then the height of $x$ on $k$ computed by $M_k$ is the same as the original height of $x$ on $k_0$ computed by $M_{k_0}$. This follows from the definition. 
\end{enumerate}

Now since $k_0$ is a product formula field which satisfies the Northcott property, we conclude that $k$ also satisfies the Nothcott property by considering the minimal polynomial of the elements and using the inequality $\mathrm{max}\{h(x+y),h(xy)\}\leq h(x)+h(y)$ of heights (on positive-characteristic fields).
\end{prf}

Next, we introduce the Weil height machinery. See \cite[Section 2.8]{Ser97} for a reference. In the following statements, we let our coefficient field $F\in\{\mathbb{R},\mathbb{C}\}$.

\begin{theorem}
Let $k$ be a product formula field and let $X$ be a projective variety over $k$. Denote $H_F$ as the quotient of the vector space of $F$-valued functions on $X(k)$ by the space of bounded $F$-valued functions on $X(k)$. Then there is a unique $F$-linear map $L\mapsto h_{L}$ of $\mathrm{Pic}(X)_{F}$ to $H_F$ such that for every morphism $\phi:X\rightarrow\mathbb{P}_{k}^{N}$, we have $h_{\phi^{*}\mathcal{O}(1)}=h_{\phi}+O(1)$ in which $h_{\phi}(x)=h(\phi(x))$ is the naive height calculated on the projective space.
\end{theorem}

Since the case where the coefficient field $F=\mathbb{C}$ is less commonly treated in the literature, we shall do some explanations. The original statement in \cite[Section 2.8]{Ser97} gives a unique group homomorphism from $\mathrm{Pic}(X)$ to $H_{\mathbb{R}}$ which satisfies the desired property. So when $F=\mathbb{R}$, we can get the $\mathbb{R}$-linear map $\mathrm{Pic}(X)_{\mathbb{R}}\rightarrow H_{\mathbb{R}}$ by tensoring $\mathbb{R}$. When $F=\mathbb{C}$, we notice that there is a natural inclusion $H_{\mathbb{R}}\hookrightarrow H_{\mathbb{C}}$ (with $H_{\mathbb{R}}$ and $H_{\mathbb{C}}$ defined as in the theorem above). So in order to get the desired $\mathbb{C}$-linear map $\mathrm{Pic}(X)_{\mathbb{C}}\rightarrow H_{\mathbb{C}}$, we shall firstly gain the group homomorphism from $\mathrm{Pic}(X)$ to $H_{\mathbb{C}}$ by $\mathrm{Pic}(X)\rightarrow H_{\mathbb{R}}\hookrightarrow H_{\mathbb{C}}$, and then tensoring $\mathbb{C}$.

The following statements are immediate from definition. We cannot state part (ii) in a unified way since the concept of ``ample $\mathbb{C}$-divisor" does not make sense. Notice that by abusing notation, we often let $h_L$ be a fixed representative of the equivalence class. We shall guarantee that there is no ambiguity.

\begin{lemma}
Let $k$ be a product formula field.
\begin{enumerate}
\item
Let $f:X\rightarrow Y$ be a morphism of projective varieties over $k$. Then $h_{f^*L}(x)=h_{L}(f(x))+O(1)$ as functions on $X(k)$ for every $L\in\mathrm{Pic}(Y)_{F}$.

\item
Let $X$ be a projective variety over $k$.

Suppose the coefficient field $F=\mathbb{R}$. Let $L$ be an ample $\mathbb{R}$-divisor on $X$. Then $h_L$ is bounded below. Suppose further that $k$ satisfies the Northcott property, then $\{x\in X(k)|\ h_{L}(x)\leq M\}$ is a finite set for every $M>0$ (and every representative of $h_{L}$).

Suppose the coefficient field $F=\mathbb{C}$. Let $L$ be an ample line bundle on $X$. Then $\mathrm{Re}(h_L)$ is bounded below. Suppose further that $k$ satisfies the Northcott property, then the set $\{x\in X(k)|\ \mathrm{Re}(h_{L}(x))\leq M\}$ is finite for every $M>0$ (and every representative of $h_{L}$).
\end{enumerate}
\end{lemma}

We shall use the fact that the height function associated to an effective line bundle is bounded below on an open dense subset. See \cite[Section 2.10]{Ser97} for a reference.

\begin{proposition}
Suppose the coefficient field $F=\mathbb{R}$. Let $k$ be a product formula field and let $X$ be a projective variety over $k$. Let $L\in\mathrm{Pic}(X)$ be an effective line bundle. Then there exists an open dense subset $U\subseteq X$ such that $h_{L}|_{U(k)}$ is bounded below.
\end{proposition}

\subsection{The Moriwaki height machinery}

In this subsection, we recall the Moriwaki height machinery \cite{Mor00}, in order to provide solid evidence for the statement in Remark 1.11. In the following, we let $k$ be a finitely generated extension field of $\mathbb{Q}$, i.e. an ``arithmetic function field".

Indeed, one can also get a height theory for projective varieties over $k$ by the Weil height machinery. But in this way the Northcott property will fail, as the ``constant field" $\mathbb{Q}$ has infinitely many elements. To overcome this problem, we need the Moriwaki heights.

We briefly review the construction of the Moriwaki height function. Let $B$ be a normal arithmetic variety, which is a model of $k$. We fix a big and nef hermitian line bundle $\overline{H}$ on $B$. Let $X$ be a projective variety over $k$ and let $H_{\mathbb{R}}$ be the same as in Theorem 2.10. Then we can define a group homomorphism from $\mathrm{Pic}(X)$ to $H_{\mathbb{R}}$ (denote as $L\mapsto h_L$) in the following way.

Let $L\in\mathrm{Pic}(X)$. Let $\pi:\mathcal{X}\rightarrow B$ be a model of $X$, in which $\mathcal{X}$ is an arithmetic variety. By choosing $\mathcal{X}$ carefully, we can find a line bundle $\mathcal{L}\in\mathrm{Pic}(\mathcal{X})$ that restricts to $L\in\mathrm{Pic}(X)$. We form $\mathcal{L}$ into a hermitian line bundle $\overline{\mathcal{L}}$ over $\mathcal{X}$. Then we define $h_L$ by the formula $h_L(x)=\overline{\mathcal{L}}\cdot\pi^*\overline{H}^d\cdot\overline{x}$ for every $x\in X(k)$, in which $\overline{x}$ is the closure of $x$ in $\mathcal{X}$ and $d$ is the transcendence degree of $k/\mathbb{Q}$. This notion is well-defined (i.e. independent of the choice of $(\mathcal{X},\overline{\mathcal{L}})$) up to bounded functions, and in this way we indeed get a group homomorphism from $\mathrm{Pic}(X)$ to $H_{\mathbb{R}}$. Notice that we have fixed $(B,\overline{H})$ in the procedure above. Since we will never modify these data, we do not emphasize them in the notion.

We remark that Moriwaki can define this notion not only for rational points, but also for algebraic points. But since we only need to use the height for rational points, we choose to focus on them in order to be parallel with subsection 2.2.

We list some properties of this map, following \cite[Proposition 3.3.7]{Mor00}. Notice that the projection formula here is a direct consequence of the projection formula in arithmetic intersection theory.

\begin{proposition}
Let $L\in\mathrm{Pic}(X)$.
\begin{enumerate}
\item
The function $h_L$ is bounded below on $(X\backslash\mathrm{Bs}(L))(k)$ in which $\mathrm{Bs}(L)$ is the base locus of $L$.

\item (Northcott)
If $L$ is ample, then the set $\{x\in X(k)|\ h_L(x)\leq M\}$ is finite for every $M>0$ (and every representative $h_L$). 

\item (projection formula)
Let $Y$ be another projective variety over $k$ and let $f:X\rightarrow Y$ be a morphism. Then for every $L\in\mathrm{Pic}(Y)$, we have $h_{f^*L}(x)=h_L(f(x))+O(1)$ as functions on $X(k)$.
\end{enumerate}
\end{proposition}

Now for $F\in\{\mathbb{R},\mathbb{C}\}$, we can get the $F$-linear map from $\mathrm{Pic}(X)_{F}$ to $H_F$ in the same way as we described in the paragraph below Theorem 2.10. Proposition 2.13 guarantees that this height machinery satisfies the properties listed in Lemma 2.11 and Proposition 2.12. In particular, it admits the Northcott property in the sense of part (ii) above. Indeed, it is these properties of the height machinery that we need to use in the rest of this article. The reader may find that the arguments in Section 3 and Section 4 remain valid almost line by line in characteristic 0, after substituting the Weil height machinery by the Moriwaki height machinery.

\section{Cohomologically hyperbolic endomorphisms}

Starting from this Section, we shall restrict ourselves in the case that the base field $K$ has characteristic $p>0$. We have explained our justification for this in Remark 1.11.

We will prove Theorem 1.1 and its corollary in this section. As we have mentioned in the Introduction, we need to find \emph{two different} speeds of growth. This is the main theme of the proof.

For a projective variety $X$, we denote $\mathrm{Pic}^0(X)\subseteq\mathrm{Pic}(X)$ as the subgroup consists of all algebraically trivial line bundles. Then there is a natural exact sequence $0\rightarrow\mathrm{Pic}^0(X)_{\mathbb{C}}\rightarrow\mathrm{Pic}(X)_{\mathbb{C}}\rightarrow\mathrm{N}^1(X)_{\mathbb{C}}\rightarrow0$ of $\mathbb{C}$-linear spaces. In addition, for any irreducible closed subcurve $C\subseteq X$, the intersection pairing gives a $\mathbb{C}$-linear map $\mathrm{Pic}(X)_{\mathbb{C}}\rightarrow\mathbb{C}$ which sends $L$ to $L\cdot C$. This map factors through $\mathrm{N}^1(X)_{\mathbb{C}}$ and hence we can talk about the intersection number $L\cdot C$ for $L\in\mathrm{Pic}(X)_{\mathbb{C}}$ or $L\in\mathrm{N}^1(X)_{\mathbb{C}}$. In particular, if $X$ is a projective curve, then we can talk about the degree of $L\in\mathrm{Pic}(X)_{\mathbb{C}}$ or $L\in\mathrm{N}^1(X)_{\mathbb{C}}$.

We begin with a proposition which shall be proved by a height argument. Whenever we use the height machinery in this section, we let the coefficient field $F=\mathbb{C}$.

\begin{proposition}
Let $X$ be a projective variety and let $f$ be an endomorphism of $X$. Let $C\subseteq X$ be an irreducible closed subcurve and let $L_1,L_2\in\mathrm{Pic}(X)_{\mathbb{C}}$ be two $\mathbb{C}$-divisors on $X$. Let $\mu_1,\mu_2\in\mathbb{C}$ be two numbers such that $|\mu_1|\neq|\mu_2|$. Suppose that
\begin{enumerate}
\item
$L_1\cdot C\neq0$ and $L_2\cdot C\neq0$, and
\item
there exists a positive integer $m$ such that $(f^*-\mu_1)^m(L_1)=(f^*-\mu_2)^m(L_2)=0$ in $\mathrm{Pic}(X)_{\mathbb{C}}$.
\end{enumerate}
Then for every point $x\in X(K)$, the set $\mathcal{O}_{f}(x)\cap C(K)$ is finite.
\end{proposition}

We need an elementary lemma about the speed of growth of differential sequences.

\begin{lemma}
Let $(x_n)_{n\geq0}$ be a sequence of complex numbers and let $a\in\mathbb{C}$. We extend the definition of $x_n$ to the negative numbers by setting $x_n=0$ for $n<0$. Then we get a sequence $(x_n)_{n\in\mathbb{Z}}$. We inductively define the differential sequences as follows.
\begin{enumerate}
\item
We let $x_n^{(0)}=x_n$ for every integer $n$.
\item
We let $x_n^{(i)}=x_n^{(i-1)}-ax_{n-1}^{(i-1)}$ for every integer $n$ and every positive integer $i$.
\end{enumerate}

Let $m$ be a positive integer such that the sequence $(x_n^{(m)})_{n\geq0}$ is bounded. Then the following holds.
\begin{enumerate}
\item
Suppose $|a|>1$ and the sequence $(x_n)_{n\geq0}$ is unbounded. Let $k$ be the maximal element in $\{0,\dots,m-1\}$ such that $(x_n^{(k)})_{n\geq0}$ is unbounded. Then the limit $\lim\limits_{n\rightarrow\infty}\frac{x_n}{n^ka^n}$ exists and is nonzero.
\item
Suppose $|a|=1$. Then there exists $C>0$ such that $|x_n|\leq Cn^m$ for every positive integer $n$.
\item
Suppose $|a|<1$. Then the sequence $(x_n)_{n\geq0}$ is bounded.
\end{enumerate}
\end{lemma}

\begin{prf}
We have the formula $x_n^{(i-1)}=x_{n}^{(i)}+ax_{n-1}^{(i)}+\cdots+a^{n-N-1}x_{N+1}^{(i)}+a^{n-N}x_{N}^{(i-1)}$ for every $n\geq N$ and every positive integer $i$. Using this formula, one can prove the assertions by induction. We will write a detailed proof for part (i), and then one can prove parts (ii)(iii) easily by using the same method.

We prove by induction on $k$. Assume $k=0$. Then the sequence $(x_n^{(1)})_{n\geq0}$ is bounded. We pick $M>0$ such that $|x_n^{(1)}|\leq M$ for every nonnegative integer $n$. Now by the formula above, we can see that $|\frac{x_n}{a^n}-\frac{x_N}{a^N}|\leq\frac{M}{|a|^n}\cdot\frac{|a|^{n-N}-1}{|a|-1}<\frac{M}{|a|-1}\cdot\frac{1}{|a|^N}$ for every $n\geq N$. Since $|a|>1$, we conclude that $(\frac{x_n}{a^n})_{n\geq0}$ is a Cauchy sequence. Hence $\lim\limits_{n\rightarrow\infty}\frac{x_n}{a^n}$ exists. Since the sequence $(x_n)_{n\geq0}$ is assumed to be unbounded, we can find a positive integer $N$ such that $|x_N|>\frac{M}{|a|-1}$. Fixing this $N$ in the formula above and let $n$ tend to infinity, we see that $\lim\limits_{n\rightarrow\infty}\frac{x_n}{a^n}\neq0$.

Now assume the assertion holds in the case $k=k_0-1$. We show that it also holds in the case $k=k_0$. Using the induction hypothesis towards the sequence $(x_n^{(1)})_{n\geq0}$, we see that the limit $\lim\limits_{n\rightarrow\infty}\frac{x_n^{(1)}}{n^{k_0-1}a^n}$ exists and is nonzero. Denote this limit by $C$. We show that $\lim\limits_{n\rightarrow\infty}\frac{x_n}{n^{k_0}a^n}=\frac{C}{k_0}$.

We denote $\delta_n=\frac{x_n^{(1)}}{n^{k_0-1}a^n}-C$ for every positive integer $n$. We fix an arbitrary $\varepsilon>0$. Let $N_0$ be a positive integer such that $|\delta_n|<\frac{\varepsilon}{2}$ for every $n\geq N_0$. By taking $N=N_0$ in the formula at the beginning, we get
\begin{equation*}
\begin{split}
x_n&=x_n^{(1)}+ax_{n-1}^{(1)}+\cdots+a^{n-N_0-1}x_{N_0+1}^{(1)}+a^{n-N_0}x_{N_0}\\&=(C+\delta_n)n^{k_0-1}a^n+\cdots+(C+\delta_{N_0+1})(N+1)^{k_0-1}a^n+\frac{x_{N_0}}{a^{N_0}}a^n\\&=\left(C\cdot(n^{k_0-1}+\cdots+(N_0+1)^{k_0-1})+(\delta_{n}n^{k_0-1}+\cdots+\delta_{N_0+1}(N_0+1)^{k_0-1})+\frac{x_{N_0}}{a^{N_0}}\right)\cdot a^n
\end{split}
\end{equation*}
for every $n\geq N_0$. Thus for every $n\geq N_0$, we have
$$
\frac{x_n}{n^{k_0}a^n}=C\cdot\frac{n^{k_0-1}+\cdots+(N_0+1)^{k_0-1}}{n^{k_0}}+\frac{\delta_{n}n^{k_0-1}+\cdots+\delta_{N_0+1}(N_0+1)^{k_0-1}}{n^{k_0}}+\frac{x_{N_0}}{n^{k_0}a^{N_0}}.
$$

Notice that the general term formula of $\sum\limits_{i=1}^{n}i^{k_0-1}$ is a polynomial of $n$ with leading term $\frac{n^{k_0}}{k_0}$. So we can conclude that there exists an integer $N\geq N_0$, such that $|\frac{x_n}{n^{k_0}a^n}-\frac{C}{k_0}|<\frac{\varepsilon}{4}+\frac{\varepsilon}{2}+\frac{\varepsilon}{4}=\varepsilon$ for every $n\geq N$. Hence we have proved that $\lim\limits_{n\rightarrow\infty}\frac{x_n}{n^{k_0}a^n}=\frac{C}{k_0}\neq0$. Thus we finish the proof by induction.
\end{prf}

Now we can prove Proposition 3.1.

\proof[Proof of Proposition 3.1]
Write $L_1=\sum\limits_{i=1}^{s}a_iM_i$ and $L_2=\sum\limits_{i=1}^{t}b_iN_i$ where $a_1,\dots,a_s,b_1,\dots,b_t\in\mathbb{C}$ and $M_1,\dots,M_s,N_1,\dots,N_t\in\mathrm{Pic}(X)$. By assumption (i), we may write $L_1|_{C}=z_1(c_0A_0+\sum\limits_{i=1}^{s'}c_iA_i)$ and $L_2|_{C}=z_2(d_0B_0+\sum\limits_{i=1}^{t'}d_iB_i)$ in which $A_0,A_1,\dots,A_{s'},B_0,B_1,\dots,B_{t'}\in\mathrm{Pic}(C)$ are very ample line bundles and $z_1,z_2,c_0,c_1,\dots,c_{s'},d_0,d_1,\dots,d_{t'}$ are complex numbers satisfying $z_1,z_2\neq0$, $\mathrm{Re}(c_0),\mathrm{Re}(d_0)>0$, and $\mathrm{Re}(c_1),\dots,\mathrm{Re}(c_{s'}),\mathrm{Re}(d_1),\dots,\mathrm{Re}(d_{t'})\geq0$. We will only explain this for $L_1|_{C}$ since the same argument also works for $L_2|_{C}$. Fix a very ample line bundle $A_0$ on $C$. Let $z_1=\frac{\text{deg}(L_1|_{C})}{\text{deg}(A_0)}\neq0$. Then $L_1|_{C}=z_1(A_0+N)$ for some $N\in\text{Pic}^0(C)_{\mathbb{C}}$. We write $N=\sum\limits_{i=1}^{s'}e_iZ_i$ for some complex numbers $e_1,\dots,e_{s'}$ and some $Z_1,\dots,Z_{s'}\in\text{Pic}^0(C)$. By changing $Z_i$ into $-Z_i$ if necessary, we may assume that the real parts of $e_1,\dots,e_{s'}$ are nonnegative. Then one can see that the assertion holds as each $Z_i$ will become ample after adding an arbitrarily small positive multiple of $A_0$.

Now suppose by contradiction that there exists a point $x_0\in X(K)$ such that $\mathcal{O}_{f}(x_0)\cap C(K)$ is an infinite set. In order to make use of the height machinery, we want to find a finitely generated field $k\subseteq K$ on which all data are defined. By the standard spreading-out argument, we can find a subfield $k$ of $K$ which is finitely generated over the prime field $\mathbb{F}_p$ such that the following holds.
\begin{enumerate}
\item
The projective variety $X$ and the endomorphism $f$ are defined over $k$.
\item
The irreducible closed subcurve $C\subseteq X$ and the starting point $x_0\in X(K)$ are defined over $k$ (as a closed subvariety of $X$ and a $k$-point in $X$).
\item
The line bundles $M_1,\dots,M_s,N_1,\dots,N_t$ and $A_0,A_1,\dots,A_{s'},B_0,B_1,\dots,B_{t'}$ introduced above are pullbacks of line bundles on the model of $X$ and $C$ respectively. Moreover, we can require that $A_0,A_1,\dots,A_{s'},B_0,B_1,\dots,B_{t'}$ are still very ample on the model of $C$.
\end{enumerate}

We regard all of the data as objects over $k$ by abusing notation as follows.
\begin{enumerate}
\item
We regard $X$ as a projective $k$-variety and $f$ as a $k$-endomorphism.
\item
We regard $C$ as an irreducible closed subcurve of $X$ (over $k$) and regard the starting point $x_0$ as an element of $X(k)$. Then $\mathcal{O}_{f}(x_0)\cap C(k)$ is still an infinite set.
\item
On the model over $k$, we still denote $L_1=\sum\limits_{i=1}^{s}a_iM_i$ and $L_2=\sum\limits_{i=1}^{t}b_iN_i$. They are elements in $\mathrm{Pic}(X)_{\mathbb{C}}$. Then the equations $(f^*-\mu_1)^m(L_1)=(f^*-\mu_2)^m(L_2)=0$, $L_1|_{C}=z_1(c_0A_0+\sum\limits_{i=1}^{s'}c_jA_j)$, and $L_2|_{C}=z_2(d_0B_0+\sum\limits_{i=1}^{t'}d_jB_j)$ still hold (in $\mathrm{Pic}(X)_{\mathbb{C}}$ and $\mathrm{Pic}(C)_{\mathbb{C}}$ respectively) because the homomorphism between Picard groups induced by the base extension is injective. Notice we have required that $A_0,A_1,\dots,A_{s'},B_0,B_1,\dots,B_{t'}\in\mathrm{Pic}(C)$ are still very ample.
\end{enumerate}

Since $x_0$ cannot be an $f$-preperiodic point in $X(k)$, we know that $k$ has a positive transcendence degree over $\mathbb{F}_p$. Using Proposition 2.9, we may equip a product formula on $k$ such that $k$ satisfies the Northcott property. Now we fix representative height functions $h_{L_1},h_{L_2}:X(k)\rightarrow\mathbb{C}$ and $h_{A_0},h_{A_1},\dots,h_{A_{s'}},h_{B_0},h_{B_1},\dots,h_{B_{t'}}:C(k)\rightarrow\mathbb{C}$. Since $A_0,A_1,\dots,A_{s'},B_0,B_1,\dots,B_{t'}$ are very ample, we may let all of the functions $h_{A_0},h_{A_1},\dots,h_{A_{s'}},h_{B_0},h_{B_1},\dots,h_{B_{t'}}$ take values on $\mathbb{R}_{\geq0}$. We can see that the functions $\sum\limits_{u=0}^{m}(-1)^u{m\choose u}\mu_1^{u}h_{L_1}(f^{m-u}(x)),\sum\limits_{u=0}^{m}(-1)^u{m\choose u}\mu_2^{u}h_{L_2}(f^{m-u}(x))$ and $h_{L_1}(x)-z_1(c_0h_{A_0}(x)+\sum\limits_{i=1}^{s'}c_ih_{A_i}(x)),h_{L_2}(x)-z_2(d_0h_{B_0}(x)+\sum\limits_{i=1}^{t'}d_ih_{B_i}(x))$ are bounded on $X(k)$ and $C(k)$ respectively by using Lemma 2.11(i).

Now we fix an ample line bundle $A$ on $C$ and fix a representative height function $h_A:C(k)\rightarrow\mathbb{C}$ which takes values on $\mathbb{R}_{\geq0}$. We say two $\mathbb{R}_{\geq0}$-valued functions $g_1$ and $g_2$ are \emph{bounded by each other} if there exist $C_1,C_2>0$ such that both $g_1\leq C_1g_2+C_2$ and $g_2\leq C_1g_1+C_2$ holds. Since $\mathrm{Re}(c_0)>0$, $\mathrm{Re}(c_1),\dots,\mathrm{Re}(c_{s'})\geq0$ and $A_0,A_1,\dots,A_{s'}$ are ample, we know that the functions $|h_A|$ and $|z_1(c_0h_{A_0}+\sum\limits_{i=1}^{s'}c_ih_{A_i})|$ are bounded by each other. For the same reason, the same result holds for the functions $|h_A|$ and $|z_2(d_0h_{B_0}+\sum\limits_{i=1}^{t'}d_ih_{B_i})|$. Therefore, we conclude that the three functions $|h_{L_1}(x)|,|h_{L_2}(x)|$ and $|h_{A}(x)|$ are bounded by each other on $C(k)$. In particular, each of them is unbounded on any infinite set since $k$ satisfies the Northcott property.

Next, we will use Lemma 3.2 towards the sequences $(h_{L_1}(f^{n}(x_0)))_{n\geq0}$ and $(h_{L_2}(f^{n}(x_0)))_{n\geq0}$. Since we have assumed $\mathcal{O}_f(x_0)\cap C(k)$ to be an infinite set, the last sentence of the paragraph above guarantees that these two sequences are unbounded. So by Lemma 3.2(iii), we get $|\mu_1|,|\mu_2|\geq1$. Since $|\mu_1|\neq|\mu_2|$, we assume that $|\mu_1|>|\mu_2|\geq1$ without loss of generality. By Lemma 3.2 again, we can pick $C_1,C_2>0$ and $N\in\mathbb{Z}_+$ such that $|h_{L_1}(f^{n}(x_0))|\geq C_1|\mu_1|^n$ and $|h_{L_2}(f^{n}(x_0))|\leq C_2n^m|\mu_2|^n$ for every integer $n\geq N$. But as $\mathcal{O}_f(x_0)\cap C(k)$ is an infinite set, this contradicts with the fact that the functions $|h_{L_1}(x)|$ and $|h_{L_2}(x)|$ are bounded by each other on $C(k)$. This contradiction finishes the proof.
\endproof

Notice that we need two $\mathbb{C}$-divisors in the hypothesis of Proposition 3.1. The next proposition is a trick which shows that we can generate the second one from the first one, at least on the level of $\mathrm{N}^1(X)_{\mathbb{C}}$.

\begin{proposition}
Let $X$ be a projective variety and let $f$ be an endomorphism of $X$. Let $C\subseteq X$ be an irreducible closed subcurve. Let $\mu\in\mathbb{C}$ be an algebraic number and let $\mu'\in\mathbb{C}$ be a Galois conjugate of $\mu$. Let $m$ be a positive integer. Suppose there exists an element $L\in\mathrm{N}^1(X)_{\mathbb{C}}$ such that $L\cdot C\neq0$ and $(f^*-\mu)^m(L)=0$ in $\mathrm{N}^1(X)_{\mathbb{C}}$. Then there exists $L'\in\mathrm{N}^1(X)_{\mathbb{C}}$ such that $L'\cdot C\neq0$ and $(f^*-\mu')^m(L')=0$ in $\mathrm{N}^1(X)_{\mathbb{C}}$.
\end{proposition}

\begin{prf}
We fix $L_1,\dots,L_{\rho}\in\mathrm{Pic}(X)$ such that they form a $\mathbb{Z}$-basis of $\mathrm{N}^1(X)$. Let $A\in M_{\rho}(\mathbb{Z})$ be the matrix of $f^*:\mathrm{N}^1(X)\rightarrow\mathrm{N}^1(X)$ with respect to this basis. The hypothesis gives us a vector $v\in\mathbb{C}^{\rho}$ such that $(L_1\cdot C,\dots,L_{\rho}\cdot C)\cdot v\neq0$ and $(\mu I_{\rho}-A)^m\cdot v=0$. Our goal is to find a vector $v'\in\mathbb{C}^{\rho}$ such that $(L_1\cdot C,\dots,L_{\rho}\cdot C)\cdot v'\neq0$ and $(\mu'I_{\rho}-A)^m\cdot v'=0$.

Firstly, we may assume that all of the coefficients of the vector $v$ above are contained in the number field $\mathbb{Q}(\mu)$. We denote $\sigma:\mathbb{Q}(\mu)\stackrel{\sim}\rightarrow\mathbb{Q}(\mu')$ as the field isomorphism which sends $\mu$ to $\mu'$. Then we have $(L_1\cdot C,\dots,L_{\rho}\cdot C)\cdot\sigma(v)\neq0$ and $(\mu'I_{\rho}-A)^{m}\cdot\sigma(v)=0$ since the numbers $L_1\cdot C,\dots,L_{\rho}\cdot C$ and the coefficients of $A$ are contained in $\mathbb{Z}$. Hence we finish the proof by taking $v'=\sigma(v)$.
\end{prf}

Proposition 3.3 reveals the property of $\mathfrak{Root}$ that we need.

\begin{lemma}
Let $\mu\in\mathbb{R}_{>0}$ be an algebraic integer. Suppose that the modulus of every Galois conjugate of $\mu$ equals to $\mu$. Then $\mu\in\mathfrak{Root}$.
\end{lemma}

\begin{prf}
Let $\mu_1=\mu,\mu_2,\dots,\mu_d$ be all of the Galois conjugates of $\mu$. Since $\mu$ is an algebraic integer, we can see that $\mu^d=|\mu_1\mu_2\cdots\mu_d|$ is a positive integer. Hence the result follows.
\end{prf}

Now, we find that the gap between $\mathrm{Pic}(X)_{\mathbb{C}}$ and $\mathrm{N}^1(X)_{\mathbb{C}}$ is a difficulty for us to prove Theorem 1.1. But notice that if the $\mathbb{C}$-linear surjection $\mathrm{Pic}(X)_{\mathbb{C}}\rightarrow\mathrm{N}^1(X)_{\mathbb{C}}$ admits an $f^*$-equivariant section, then the difference between $\mathrm{Pic}(X)_{\mathbb{C}}$ and $\mathrm{N}^1(X)_{\mathbb{C}}$ can be eliminated for our purpose. The next lemma shows that such a section exists in certain cases.

\begin{lemma}
Let $X$ be a projective variety and let $f$ be an endomorphism of $X$. Then the $\mathbb{C}$-linear surjection $\mathrm{Pic}(X)_{\mathbb{C}}\rightarrow\mathrm{N}^1(X)_{\mathbb{C}}$ admits an $f^*$-equivariant section in the following two cases.
\begin{enumerate}
\item
$\mathrm{Pic}^0(X)=\{0\}$.
\item
$X$ is an abelian variety and $f$ is a group endomorphism of $X$.
\end{enumerate}
\end{lemma}

\begin{prf}
In case (i), the surjection $\mathrm{Pic}(X)_{\mathbb{C}}\rightarrow\mathrm{N}^1(X)_{\mathbb{C}}$ is indeed an isomorphism of $\mathbb{C}$-linear spaces and hence the assertion holds. So we may restrict ourselves to the setting of case (ii) from now on. We need the following facts about $\mathrm{Pic}^0(X)$, which can be learned from \cite[Section 8]{Mum08}.
\begin{enumerate}
\item
For any $L\in\mathrm{Pic}(X)$, we have $L-[-1]^*L\in\mathrm{Pic}^0(X)$.
\item
For any $L\in\mathrm{Pic}^0(X)$, we have $L=-[-1]^*L$.
\end{enumerate}

For $L\in\mathrm{Pic}(X)_{\mathbb{C}}$, we say that $L$ is \emph{symmetric} if $L=[-1]^*L$. By the two facts above, one can prove that there is exactly one symmetric element in each fiber of the surjection $\mathrm{Pic}(X)_{\mathbb{C}}\rightarrow\mathrm{N}^1(X)_{\mathbb{C}}$. Indeed, for any $L\in\mathrm{Pic}(X)_{\mathbb{C}}$, the $\mathbb{C}$-divisor $\frac{1}{2}(L+[-1]^*L)$ is symmetric and lies in the same fiber with $L$ by fact (i). On the other hand, the uniqueness is guaranteed by fact (ii) because it implies that there is no nonzero symmetric element in $\mathrm{Pic}^0(X)_{\mathbb{C}}$.

Now we consider the section $\mathrm{N}^1(X)_{\mathbb{C}}\rightarrow\mathrm{Pic}(X)_{\mathbb{C}}$ which sends each element to the unique symmetric element in the corresponding fiber. Then this section is indeed $\mathbb{C}$-linear, and it is also $f^*$-equivariant because $f$ is a group endomorphism of $X$. Thus we finish the proof.
\end{prf}

\begin{remark}
The careful reader may find that the definition of $\mathrm{Pic}^0(X)$ for abelian varieties in \cite[Section 8]{Mum08} is seemingly different to ours. But according to part (vi) on page 71 and Theorem 1 on page 73 in loc. cit., one can see that the two definitions are equivalent.
\end{remark}

We can prove some special cases of Theorem 1.1 now. The general case will be proved as a consequence of these two special cases by using an Albanese argument.

\begin{proposition}
Let $X$ be a projective variety and let $f$ be a surjective endomorphism of $X$. Let $C\subseteq X$ be an irreducible closed subcurve and let $x\in X(K)$ be a point. Suppose that $\mu_i(f)\notin\mathfrak{Root}$ for every $i\in\{1,\dots,\mathrm{dim}(X)\}$.
\begin{enumerate}
\item
If $\mathrm{Pic}^0(X)=\{0\}$ and $\overline{\mathcal{O}_f(x)}=X$, then the set $\mathcal{O}_{f}(x)\cap C(K)$ is finite.
\item
If $X$ is an abelian variety and $f$ is a group endomorphism of $X$, then the set $\mathcal{O}_{f}(x)\cap C(K)$ is finite. We do \emph{not} need to assume $\overline{\mathcal{O}_f(x)}=X$ in this case.
\end{enumerate}
\end{proposition}

\begin{prf}
\begin{enumerate}
\item
Using Theorem 2.3, we can find $L_1,\dots,L_{\text{dim}(X)}\in\text{N}^1(X)_{\mathbb{R}}$ such that $L_1+\cdots+L_{\text{dim}(X)}\in\text{Big}(X)$ and $(f^*-\mu_1(f))^{\rho}(L_1)=\cdots=(f^*-\mu_{\text{dim}(X)}(f))^{\rho}(L_{\text{dim}(X)})=0$ in $\text{N}^1(X)_{\mathbb{R}}$. Then there exists a proper closed subset $E\subseteq X$ such that $(L_1+\cdots+L_{\text{dim}(X)})\cdot C'>0$ for every irreducible closed subcurve $C'\subseteq X$ which is not contained in $E$.

Suppose by contradiction that $\mathcal{O}_{f}(x)\cap C(K)$ is an infinite set. After substituting $C$ by an appropriate iteration $f^{n_0}(C)$, we can assume that $C$ is not contained in $E$ because $\overline{\mathcal{O}_f(x)}=X$. Then there exists $L\in\text{N}^1(X)_{\mathbb{R}}$ and $\mu\in\{\mu_1(f),\dots,\mu_{\text{dim}(X)}(f)\}$ such that $L\cdot C\neq 0$ and $(f^*-\mu)^{\rho}(L)=0$ in $\text{N}^1(X)_{\mathbb{R}}$. By Theorem 2.2, we see that $\mu\in\mathbb{R}_{>0}$ is an algebraic integer. So by Lemma 3.4 and the hypothesis that $\mu\notin\mathfrak{Root}$, we can pick a Galois conjugate $\mu'\in\mathbb{C}$ of $\mu$ such that $|\mu|\neq|\mu'|$. We regard $L$ as an element in $\text{N}^1(X)_{\mathbb{C}}$. By using Proposition 3.3, we can find $L'\in\text{N}^1(X)_{\mathbb{C}}$ such that $L'\cdot C\neq0$ and $(f^*-\mu')^{\rho}(L')=0$ in $\text{N}^1(X)_{\mathbb{C}}$. Then we can lift $L$ and $L'$ into $\text{Pic}(X)_{\mathbb{C}}$ by using Lemma 3.5, and finally deduce a contradiction by Proposition 3.1. Hence we finish the proof.

\item
The proof in this case is same as above. We will only explain why we do not need to assume $\overline{\mathcal{O}_f(x)}=X$ here. In the proof above, we need this hypothesis because we may need to alter the curve $C$. But since the big cone and the ample cone are the same for abelian varieties, we do not need to do such alterations in this case. So we do not need to assume $\overline{\mathcal{O}_f(x)}=X$.
\end{enumerate}
\end{prf}

Now we deduce Theorem 1.1 from Proposition 3.7.

\proof[Proof of Theorem 1.1]
We may assume that $C$ is irreducible, and we may also assume that $X$ is normal by taking normalization. Then by Theorem 2.6, we may further assume that $f$ admits a fixed point by iterating $f$. Notice that Proposition 2.1(ii)(iv) guarantees that these procedures will not affect the hypothesis $\{\mu_1(f),\dots,\mu_{\text{dim}(X)}(f)\}\cap\mathfrak{Root}=\emptyset$.

Let $x_0\in X(K)$ be a fixed point of $f$. Then $f$ is an endomorphism of the pointed normal projective variety $(X,x_0)$. We shall consider the Albanese map $\phi:X\rightarrow A$ with respect to the point $x_0$ and we use the Appendix of $\cite{Moc12}$ as a reference for the general facts about the Albanese map. If the map $\phi$ is constant, then we know $A=0$ by the universal property and hence we have $\mathrm{Pic}^0(X)=\{0\}$ by $\cite[\mathrm{Proposition}\ \mathrm{A}.6]{Moc12}$. Thus we finish the proof by Proposition 3.7(i). So we may assume that $\phi$ is non-constant without loss of generality.

By the universal property, we can see that $f$ induces a group endomorphism $g:A\rightarrow A$ which satisfies $g\circ\phi=\phi\circ f$. We prove that $g$ is surjective and satisfies $\{\mu_1(g),\dots,\mu_{\text{dim}(A)}(g)\}\cap\mathfrak{Root}=\emptyset$. Indeed, we know that there exists a positive integer $N$ such that the map $\zeta_N:X^N\rightarrow A$ given by $(x_1,\dots,x_N)\mapsto\sum\limits_{i=1}^{N}\phi(v_j)$ is surjective by $\cite[\mathrm{Proposition}\ \mathrm{A}.3\mathrm{(ix)}]{Moc12}$. So in view of Proposition 2.1(iii)(v), we conclude that $g$ is surjective and satisfies $\{\mu_1(g),\dots,\mu_{\text{dim}(A)}(g)\}\cap\mathfrak{Root}=\emptyset$ since $g\circ\zeta_{N}=\zeta_{N}\circ(f\times\cdots\times f)$ in which $f\times\cdots\times f$ is the split endomorphism of $X^N$ induced by $f$.

Now by contradiction, we assume that $\mathcal{O}_{f}(x)\cap C(K)$ is an infinite set. Then the return set $\{n\in\mathbb{N}|\ g^{n}(\phi(x))\in\phi(C)(K)\}$ is infinite. Since $\phi(C)\subseteq A$ is either an irreducible closed subcurve or a closed point, we conclude by Proposition 3.7(ii) that $\phi(x)$ is $g$-preperiodic. But $\overline{\mathcal{O}_f(x)}=X$ implies that $\phi(X)\subseteq\overline{\mathcal{O}_g(\phi(x))}$. So we get a contradiction as we have assumed that $\phi$ is non-constant. This contradiction finishes the proof.
\endproof

Now we prove the corollary.

\proof[Proof of Corollary 1.3]
We recall some facts before the proof.
\begin{enumerate}[label=(\alph*)]
\item
We have $\mu_1(f)\mu_2(f)\cdots\mu_{\text{dim}(X)}(f)=\lambda_{\text{dim}(X)}(f)=\text{deg}(f)\in\mathbb{Z}_+$.
\item
The multipliers $\mu_1(f),\dots,\mu_{\text{dim}(X)}(f)$ are algebraic integers. See Theorem 2.2.
\item
Let $a\in\mathfrak{Root}$ and let $n$ be a positive integer. Suppose that $\frac{n}{a}$ is an algebraic integer. Then $\frac{n}{a}\in\mathfrak{Root}$. In particular, we have $a\leq n$. This is a direct consequence of the fact that a rational algebraic integer is an integer.
\end{enumerate}

Now we start the proof.

\begin{enumerate}
\item
Since $\text{dim}(X)=2$, we can assume $\overline{\mathcal{O}_f(x)}=X$ in the procedure of proving $(X,f)$ satisfies the $\text{DML}_0$ property. Using the assumption $\mu_1(f)=\lambda_1(f)\notin\mathfrak{Root}$ and the facts above, we can see that $\mu_2(f)\notin\mathfrak{Root}$ as well. Hence the result follows from Theorem 1.1.

\item
By Theorem 1.1, we only need to verify that $\{\mu_1(f),\dots,\mu_{\text{dim}(X)}(f)\}\cap\mathfrak{Root}=\emptyset$. Assume the contrary. Then there exists a multiplier $\mu_k(f)\in\mathfrak{Root}$, for which we have $\mu_k(f)\geq1$. But by the facts above, we can see that $\mu_k(f)\leq\text{deg}(f)$ as well. So we get a contradiction and thus finish the proof.

\item
Since the iterations and the subsystems of amplified endomorphisms are still amplified, we may assume that $\overline{\mathcal{O}_f(x)}=X$ by applying standard arguments. Then the assertion follows from Theorem 2.5 and part (ii).
\end{enumerate}
\endproof

\section{Split endomorphisms}

In this section, we will prove Theorem 1.8 and Corollary 1.9. It is easy to see that there are two different speeds of growth in such settings. For example, the growth on the component $C$ is quicker than that on the component $X$ in the setting of Theorem 1.8. We will use this observation to give the proofs using a height argument. We notice that whenever we use the height machinery in this section, we let the coefficient field $F=\mathbb{R}$.

We shall firstly prove the more generalized Proposition 4.1. Then we prove Theorem 1.8 and Corollary 1.9 as corollaries of Proposition 4.1. To clarify the structure, we remark that the Lemmas 4.2 and 4.3 are contained inside the proof of Proposition 4.1.

\begin{proposition}
Let $X$ and $Y_1,Y_2$ be projective varieties. Let $f:X\rightarrow X$ and $g_1:Y_1\rightarrow Y_1,g_2:Y_2\rightarrow Y_2$ be surjective endomorphisms. Let $p_1:X\rightarrow Y_1$ and $p_2:X\rightarrow Y_2$ be morphisms such that $p_1\circ f=g_1\circ p_1$ and $p_2\circ f=g_2\circ p_2$. Suppose there exists an ample $\mathbb{R}$-divisor $L_2\in\mathrm{Pic}(Y_2)_{\mathbb{R}}$ such that $g_2^*L_2-\lambda_1(g_1)L_2$ is ample. Let $V\subseteq X$ be an irreducible closed subvariety such that $\mathrm{dim}(V)=\mathrm{dim}(p_1(V))\geq1$ and let $x\in X(K)$ be a point such that $\overline{\mathcal{O}_{f}(x)\cap V}=V$. Then $p_2(x)$ is a $g_2$-preperiodic point and hence $V\subseteq p_2^{-1}(y_2)$ for some point $y_2\in Y_2(K)$.
\end{proposition}

\begin{prf}
Since the morphism from $V\subseteq X$ to $Y_1$ is supposed to be generically finite onto its image, we can find an ample line bundle $L_1\in\mathrm{Pic}(Y_1)$ such that there exist an ample line bundle $L\in\mathrm{Pic}(V)$ and an effective line bundle $E\in\mathrm{Pic}(V)$ which satisfy $(p_1^{*}L_1)|_{V}=L+E$. Moreover, we can write $L_2=\sum\limits_{i=1}^{m}a_iA_i$ and $g_2^{*}L_2-\lambda_1(g_1)L_2=\sum\limits_{j=1}^{n}b_jB_j$ for some $a_1,\dots,a_m,b_1,\dots,b_n>0$ and some ample line bundles $A_1,\dots,A_m,B_1,\dots,B_n\in\mathrm{Pic}(Y_2)$. In order to make use of the Weil height machinery, we have to find a finitely generated field $k\subseteq K$ on which all data are defined. By the standard spreading-out argument, we can find a subfield $k$ of $K$ which is finitely generated over $\mathbb{F}_p$ such that the following holds.

\begin{enumerate}
\item
The projective varieties $X,Y_1,Y_2$ and the morphisms $f,g_1,g_2,p_1,p_2$ are defined over $k$.
\item
The irreducible closed subvariety $V\subseteq X$ and the starting point $x\in X(K)$ are defined over $k$ (as a closed subvariety of $X$ and a $k$-point in $X$).
\item
The line bundles $L_1,L,E,A_1,\dots,A_m,B_1,\dots,B_n$ introduced above are pullbacks of line bundles on the model of $Y_1,V$ and $Y_2$ respectively.
\end{enumerate}

We regard all of the data as objects over $k$ by abusing notation as follows.
\begin{enumerate}
\item
We regard $X,Y_1,Y_2$ as projective $k$-varieties and $f,g_1,g_2,p_1,p_2$ as $k$-morphisms. Then $f,g_1,g_2$ are still surjective and the equations $p_1\circ f=g_1\circ p_1$ and $p_2\circ f=g_2\circ p_2$ still hold.
\item
We regard $V$ as an irreducible closed subvariety of $X$ (over $k$) and regard the starting point $x$ as an element of $X(k)$. Then $V$ is still of positive dimension and $\overline{\mathcal{O}_{f}(x)\cap V}=V$ still holds.
\item
The line bundles $L_1,L,A_1,\dots,A_m,B_1,\dots,B_n$ on $Y_1,V$ and $Y_2$ are still ample and $E\in\mathrm{Pic}(V)$ is still effective. Let $L_2=\sum\limits_{i=1}^{m}a_iA_i$ and $L_2'=\sum\limits_{j=1}^{n}b_jB_j$ be ample $\mathbb{R}$-divisors in $\mathrm{Pic}(Y_2)_{\mathbb{R}}$. Then the equations $(p_1^{*}L_1)|_{V}=L+E$ and $g_2^{*}L_2-\lambda_1(g_1)L_2=L_2'$ still hold (in $\mathrm{Pic}(V)$ and $\mathrm{Pic}(Y_2)_{\mathbb{R}}$ respectively) because the homomorphisms between Picard groups induced by the base extension are injective.
\end{enumerate}

Our goal is to prove that $p_2(x)\in Y_2(k)$ is $g_2$-preperiodic. 

Since $\text{dim}(V)\geq1$ and $\overline{\mathcal{O}_{f}(x)\cap V}=V$, the point $x\in X(k)$ cannot be $f$-preperiodic. Thus $k$ has a positive transcendence degree over $\mathbb{F}_p$. Using Proposition 2.9, we may equip a product formula on $k$ such that $k$ satisfies the Northcott property. The philosophy of our proof is quite easy: if $p_2(x)\in Y_2(k)$ is not $g_2$-preperiodic, then the growth of $\{h_{L_1}(g_1^{n}(p_1(x)))\}_{n\in\mathbb{N}}$ is slower than the growth of $\{h_{L_2}(g_2^{n}(p_2(x)))\}_{n\in\mathbb{N}}$. But by looking at an infinite subsequence of them on $V$, we find that the former can control the latter and hence get a contradiction.

To realize this philosophy, we investigate the growth of the two sequences above. Since both $L_1\in\mathrm{Pic}(Y_1)$ and $L_2\in\mathrm{Pic}(Y_2)_{\mathbb{R}}$ are ample, we may fix representatives $h_{L_1}$ and $h_{L_2}$ such that they take values in $\mathbb{R}_{\geq1}$. Firstly, we show that $\limsup\limits_{n\rightarrow\infty}h_{L_1}(g_1^{n}(p_1(x)))^{\frac{1}{n}}\leq\lambda_1(g_1)$. This is known as the Kawaguchi--Silverman--Matsuzawa's upper bound (see \cite{KS16,Mat20}) and a proof in arbitrary characteristic can be found in $\cite[\text{Proposition}\ 2.10]{Xie23}$. We include a proof here for completeness as the proof is rather easy in the case of endomorphisms of projective varieties.

\begin{lemma}
We have $\limsup\limits_{n\rightarrow\infty}h_{L_1}(g_1^{n}(p_1(x)))^{\frac{1}{n}}\leq\lambda_1(g_1)$.
\end{lemma}

\begin{prf}
We will prove that $\limsup\limits_{n\rightarrow\infty}h_{L_1}(g_1^{n}(p_1(x)))^{\frac{1}{n}}\leq\lambda_1(g_1)+\varepsilon$ for every $\varepsilon>0$. Recall that $\lambda_1(g_1)$ is the spectral radius of the action $g_1^*$ on $\mathrm{N}^{1}(Y_{1,K})_{\mathbb{R}}$. So we may find a sequence $\{a_m\}_{m\in\mathbb{N}}$ of integers $\geq2$ and a number $C>0$ such that
\begin{enumerate}
\item
$a_m\leq C\cdot(\lambda_1(g_1)+\varepsilon)^m$ for any $m\in\mathbb{N}$, and
\item
$a_mL_1-(g_1^m)^*L_1$ is an ample line bundle on $Y_1$ for any $m\in\mathbb{N}$
\end{enumerate}
as $L_1\in\mathrm{Pic}(Y_1)$ is ample. In order to prove that $\limsup\limits_{n\rightarrow\infty}h_{L_1}(g_1^{n}(p_1(x)))^{\frac{1}{n}}\leq\lambda_1(g_1)+\varepsilon$, we only need to show that $\limsup\limits_{n\rightarrow\infty}h_{L_1}(g_1^{n}(p_1(x)))^{\frac{1}{n}}\leq a_m^{\frac{1}{m}}$ for every $m\in\mathbb{Z}_{+}$.

Now fix $m\in\mathbb{Z}_+$. By the properties of the height machinery, we can find $C'>0$ such that $a_mh_{L_1}(y)-h_{L_1}(g_1^m(y))\geq-C'$ for every $y\in Y_1(k)$. So we have $h_{L_1}(g_1^{mN}(y))\leq a_m^{N}(h_{L_1}(y)+C')$ for every $N\in\mathbb{Z}_+$ and every $y\in Y_1(k)$. Hence $\limsup\limits_{n\rightarrow\infty}h_{L_1}(g_1^{n}(p_1(x)))^{\frac{1}{n}}\leq a_m^{\frac{1}{m}}$ and thus we finish the proof.
\end{prf}

Next, we consider the growth of $\{h_{L_2}(g_2^{n}(p_2(x)))\}_{n\in\mathbb{N}}$.

\begin{lemma}
Suppose $p_2(x)\in Y_2(k)$ is not $g_2$-preperiodic. Then there exists $C_0,\varepsilon_0>0$ such that $h_{L_2}(g_2^{n}(p_2(x)))\geq C_0(\lambda_1(g_1)+\varepsilon_0)^n$ for every $n\in\mathbb{N}$.
\end{lemma}

\begin{prf}
Since $g_2^*L_2-\lambda_1(g_1)L_2$ is an ample $\mathbb{R}$-divisor on $Y_2$, we know that there exists $\varepsilon_0>0$ such that $g_2^*L_2-(\lambda_1(g_1)+2\varepsilon_0)L_2$ is also an ample $\mathbb{R}$-divisor. So the function $h_{L_2}(g_2(y))-(\lambda_1(g_1)+2\varepsilon_0)h_{L_2}(y)$ is bounded below on $Y_2(k)$. We pick $C>0$ such that $h_{L_2}(g_2(y))\geq(\lambda_1(g_1)+2\varepsilon_0)h_{L_2}(y)-C$ for every $y\in Y_2(k)$. Since $L_2$ is an ample $\mathbb{R}$-divisor on $Y_2$, the point $p_2(x)\in Y_2(k)$ is not $g_2$-preperiodic, and $k$ satisfies the Northcott property, we can find $n_0\in\mathbb{Z}_+$ such that $h_{L_2}(g_2^{n_0}(p_2(x)))\geq\frac{C}{\varepsilon_0}$ by Lemma 2.11(ii). Then we conclude that $h_{L_2}(g_2^{n+n_0}(p_2(x)))\geq(\lambda_1(g_1)+\varepsilon_0)^nh_{L_2}(g_2^{n_0}(p_2(x)))$ for every $n\in\mathbb{N}$ and hence finish the proof.
\end{prf}

Now we can finish the proof of Proposition 4.1. Assume by contradiction that $p_2(x)\in Y_2(k)$ is not $g_2$-preperiodic.

Recall that $(p_1^{*}L_1)|_{V}=L+E$ for some ample line bundle $L\in\mathrm{Pic}(V)$ and effective line bundle $E\in\mathrm{Pic}(V)$. We fix a positive integer $M$ such that $ML-(p_2^{*}L_2)|_{V}\in\mathrm{Pic}(V)_{\mathbb{R}}$ is ample. Let $U\subseteq V$ be the open dense subset corresponding to the effective line bundle $ME$ in Proposition 2.12, then we can see that the function $Mh_{L_1}(p_1(v))-h_{L_2}(p_2(v))$ is bounded below on $U(k)$. Since $\overline{\mathcal{O}_f(x)\cap V}=V$ and $V$ is a positive dimensional irreducible $k$-variety, we know that $\{n\in\mathbb{N}|\ f^{n}(x)\in U(k)\}$ is an infinite set. But this leads to a contradiction by taking Lemmas 4.2 and 4.3 into account. So we have proved that $p_2(x)\in Y_2(k)$ must be $g_2$-preperiodic. Going back to the level of $K$, we know that $p_2(x)$ is also $g_2$-preperiodic.

The assertion about $V$ is then an immediate consequence.
\end{prf}

\begin{remark}
\begin{enumerate}
\item
By looking into the proof of \cite[Theorem 1.1]{Meng20}, one can see that our hypothesis of $g_2$ in Proposition 4.1 is equivalent to saying that the modulus of every eigenvalue of $g_2^*:\mathbb{N}^1(Y_2)_{\mathbb{R}}\rightarrow\mathbb{N}^1(Y_2)_{\mathbb{R}}$ is greater than $\lambda_1(g_1)$.

\item
Assume that the Kawaguchi--Silverman conjecture \cite[Conjecture 6]{KS16} holds for $g_2$. Then we can weaken the condition of $g_2$ into $\lambda_1(g_2)>\lambda_1(g_1)$, and the conclusion turns into saying that $\mathcal{O}_{g_2}(p_2(x))$ is not dense in $Y_2$. In view of \cite[Theorem 1.5]{Xieb}, this is indeed a strengthening of our result.
\end{enumerate}
\end{remark}

Now we prove Theorem 1.8.

\proof[Proof of Theorem 1.8]
In order to prove that $(X\times C,f\times g)$ satisfies the $\text{DML}_{\epsilon}$ property, we only need to show that for every $(x,y)\in(X\times C)(K)$ and every irreducible closed subvariety $V\subseteq X\times C$ of positive dimension, the set $\{n\in\mathbb{N}|\ (f^{n}(x),g^{n}(y))\in V(K)\}$ is of type $\epsilon$ if $\overline{\mathcal{O}_{f\times g}((x,y))\cap V}=V$. Notice that the condition $\lambda_1(f)<\mathrm{deg}(g)$ ensures that $g^{*}L-\lambda_1(f)L$ is an ample $\mathbb{R}$-divisor on $C$ for every ample line bundle $L\in\mathrm{Pic}(C)$. Thus we may apply Proposition 4.1. Let $\mathrm{pr}_{X}:X\times C\rightarrow X$ and $\mathrm{pr}_{C}:X\times C\rightarrow C$ be the two projections.
\begin{enumerate}
\item
Suppose $\mathrm{dim}(\mathrm{pr}_{X}(V))<\mathrm{dim}(V)$. Then we have $V=\mathrm{pr}_{X}(V)\times C$ as a closed subvariety of $X\times C$. Therefore, we finish the proof as we have assumed that $(X,f)$ satisfies the $\text{DML}_{\epsilon}$ property.
\item
Suppose $\mathrm{dim}(\mathrm{pr}_{X}(V))=\mathrm{dim}(V)$. Then Proposition 4.1 says that $V\subseteq\mathrm{pr}_{C}^{-1}(c)$ for some $c\in C(K)$. So $V\subseteq X\times C$ has the form $V_0\times\{c\}$ for some closed subvariety $V_0\subseteq X$. Thus we also finish the proof because $(X,f)$ satisfies the $\text{DML}_{\epsilon}$ property and the intersection of two type $\epsilon$ sets still has type $\epsilon$.
\end{enumerate}

Combining the two cases and then we are done.
\endproof

Now we turn to the proof of Corollary 1.9. The proof is by induction on $n$ and the induction step is done by Theorem 1.8. So we only need to focus on the inductive foundation.

\proof[Proof of Corollary 1.9]
By taking normalization, we may assume that $C_1,\dots,C_n$ are smooth projective curves, and then $g_1,\dots,g_n$ become surjective endomorphisms. By induction, we only need to deal with the case in which $m=n$ (i.e. all of $g_1,\dots,g_n$ are automorphisms) in view of Theorem 1.8. More concretely, for every $k\in\{m,\dots,n-1\}$, we shall apply Theorem 1.8 towards $X=C_1\times\cdots\times C_k$ and $C=C_{k+1}$. This is because we know $\lambda_1(g_1\times\cdots\times g_k)=\mathrm{deg}(g_k)<\mathrm{deg}(g_{k+1})$ according to Proposition 2.1(iii). So we reduce to the case that $m=n$.

Now the ``moreover" part is true as the 0-DML conjecture for automorphisms have been proved in \cite[Theorem 1.3]{BGT10}. Also, the assertion about the $\text{DML}_2$ property is a special case of \cite[Theorem 1.5]{XY}.
\endproof

At the end of this section, we include a proposition which says that in some cases this dynamical system satisfies the $\text{DML}_1$ property.

\begin{proposition}
Suppose $\mathrm{char}(K)=p>0$. Let $C_1,\dots,C_n,g_1,\dots,g_n$ be as in Corollary 1.9. Let $V\subseteq C_1\times\cdots\times C_n$ be a closed subvariety and let $\widetilde{C_1},\dots,\widetilde{C_m}$ be the normalization of $C_1,\dots,C_m$, respectively. Suppose one of the following conditions holds.
\begin{enumerate}
\item
None of the genus $g(\widetilde{C_i})(1\leq i\leq m)$ equals to 0.
\item
None of the genus $g(\widetilde{C_i})(1\leq i\leq m)$ equals to 1.
\item
$\mathrm{dim}(V)\leq2$.
\end{enumerate}
Then every return set of a well-defined orbit in $(C_1\times\cdots\times C_n,g_1\times\cdots\times g_n)$ into $V$ is a $p$-normal set in $\mathbb{N}$.
\end{proposition}

\begin{prf}
Suppose either (i) or (ii) holds. Then our goal is to prove that the dynamical system $(C_1\times\cdots\times C_n,g_1\times\cdots\times g_n)$ satisfies the $\text{DML}_1$ property. As above, we may assume that $C_1,\dots,C_n$ are smooth and $g_1,\dots,g_n$ are endomorphisms and thus reduce to the case in which $m=n$ by Theorem 1.8. Since the automorphism group of a smooth projective curve whose genus greater than 1 is finite, we may assume that either all of $C_1,\dots,C_m$ have genus 0 or all of them have genus 1 according to the hypothesis.

If all of $C_1,\dots,C_m$ have genus 0, then we can see that $(C_1\times\cdots\times C_m,g_1\times\cdots\times g_m)$ has the $\text{DML}_1$ property by $\cite[\text{Theorem}\ 1.8]{Der07}$. Another way to view it is that such a split automorphism of $\mathbb{P}_{K}^1\times\cdots\times\mathbb{P}_{K}^1$ is given by an \emph{affine} group action and then we reduce to the case of translation of tori by the arguments in $\cite{XY}$. If all of $C_1,\dots,C_m$ have genus 1, then they are elliptic curves and the automorphism $g_1\times\cdots\times g_m$ becomes a translation of the abelian variety $C_1\times\cdots\times C_m$ after a certain time of iterate. Thus we conclude that $(C_1\times\cdots\times C_m,g_1\times\cdots\times g_m)$ has the $\text{DML}_1$ property by \cite[Remark 3.4]{XY}.

The case in which (iii) holds is an immediate consequence of Remark 1.10.
\end{prf}

\begin{remark}
If none of the three hypotheses in Proposition 4.5 holds, then \cite[Example 5.4]{XY} reveals that the conclusion may be false.
\end{remark}

\section{Examples on the dark side}

The authors once thought that the $p$-sets in the $p$DML problem come from bounded-degree dynamical systems in some sense. So we investigated the $p$DML problem for bounded-degree systems in \cite{XY}. However, the examples in this section show that we were too naive.

Recall that in order to make the height argument work, we need to find two different speeds of growth. Generally speaking, this goal is hard to fulfill in the following cases.
\begin{enumerate}
\item
The dynamical system is a certain int-amplified endomorphism. This concept means that $f$ is an endomorphism of a projective variety $X$, for which there exists an \emph{ample} line bundle $L$ on $X$ such that $f^*L-L$ is also ample. See \cite{Meng20}. A philosophy in \cite{Xieb} says that such maps should be the algebraic analogy of expanding maps. Although some special cases can be covered by results in the previous sections, we cannot deal with the case in which the dynamical system is, for example, polarized.

\item
The dynamical system has zero entropy. In other words, the first dynamical degree $\lambda_1(f)=1$.
\end{enumerate}

We remark that ``bounded-degree" is a much stronger requirement than ``with zero entropy", and it seems that the int-amplified endomorphisms have nothing to do with bounded-degree systems. In these two cases, the height argument is hard to be applied as the dynamical system is somehow ``isotropic". But in some case, a similar height argument can work due to the additional structures of the dynamical system. Indeed, we shall deal with the automorphisms of projective surfaces in \cite{XYb}. We will treat a zero entropy case in there by using a height argument.

In this section, we proceed as follows. In subsection 5.1, we give examples of int-amplified dynamical systems which do not satisfy the $\text{DML}_0$ property. We can even construct a polarized endomorphism of $\mathbb{P}^2$ for which the $\text{DML}_0$ property fails. These examples are based on the Frobenius endomorphism in positive characteristic. Then in subsection 5.2, we give examples of dynamical systems with zero entropy for which the return sets may have a formidable form.

\subsection{Examples related to the Frobenius}

The Frobenius endomorphism can come into the picture as the following example shows.

\begin{example}
\begin{enumerate}
\item
Let the base field $K=\overline{\mathbb{F}_{p}(t)}$. Let $f:\mathbb{A}^3\rightarrow\mathbb{A}^3$ be the endomorphism given by $(x,y,z)\mapsto(x^p,(x+1)^py^p,x^pz^p)$. Let the starting point $\alpha=(t,1,1)$ and the closed subvariety $V\subseteq\mathbb{A}^3$ be the hyperplane $y=z+1$. Then we may calculate that $f^{n}(\alpha)=(t^{p^n},(t+1)^{np^n},t^{np^n})$ for all $n\in\mathbb{N}$ and hence $\{n\in\mathbb{N}|\ f^{n}(\alpha)\in V(K)\}=\{p^m|\ m\in\mathbb{N}\}$ is a ``$p$-set".

\item
A more theoretical explanation is as follows. Suppose $(X,f)$ is an isotrivial dynamical system, i.e. both the variety $X$ and the endomorphism $f:X\rightarrow X$ are defined over a finite field $\mathbb{F}_q$. Then $f$ commutes with $\mathrm{Frob}_q:X\rightarrow X$ and we let $g=\mathrm{Frob}_{q}\circ f$. Suppose further that the closed subvariety $V\subseteq X$ is also defined over $\mathbb{F}_q$. Then for any $x\in X(K)$, we have that $\{n\in\mathbb{N}|\ f^{n}(x)\in V(K)\}=\{n\in\mathbb{N}|\ g^{n}(x)\in V(K)\}$.
\end{enumerate}
\end{example}

The machinery in part (ii) above can construct int-amplified dynamical systems which do not satisfy the $\text{DML}_0$ property. For giving a concrete example, we firstly recall the following example constructed in \cite[Example 3.6]{XY}.

\begin{example}
Let $p=5$ and let $K=\overline{\mathbb{F}_{p}(t)}$. Let $E$ be the elliptic curve $x_{1}^{2}x_{2}=x_{0}^{3}+x_{2}^{3}$ in $\mathbb{P}_{K}^{2}$ with zero point $O=[0,1,0]\in E(K)$. Let $A=E\times E$ be an abelian variety. We embed $A$ into $\mathbb{P}_{K}^{8}$ by Segre embedding, i.e. $[x_0,x_1,x_2]\times[y_0,y_1,y_2]\mapsto[x_0y_0,x_0y_1,x_0y_2,x_1y_0,x_1y_1,x_1y_2,x_2y_0,$
$x_2y_1,x_2y_2]$. Let $z_{ij}$ be the coordinate of $\mathbb{P}^{8}$ corresponding to $x_iy_j$ for any $0\leq i,j\leq2$. Let $V\subseteq A$ be the closed subvariety $\{z_{02}=z_{20}+z_{22}\}\cap A$. Then there exists a point $P\in A(K)$ such that the set $\{n\in\mathbb{N}|\ n\cdot P\in V(K)\}$ is \emph{not} a finite union of arithmetic progressions in $\mathbb{N}$.
\end{example}

Now we can propose our example.

\begin{example}
Let $p,K$ and $A,V,P$ be as in Example 5.2. Let $X=A\times A$ and let $f:X\rightarrow X$ be the automorphism given by the formula $(a,b)\mapsto(a,a+b)$. Let $Y=A\times V$ be a closed subvariety of $X$ and let the starting point $x=(P,0)\in X(K)$. We can see that except for the starting point $x$, all of the data above are defined over $\mathbb{F}_p$. As a result, we may let $g=\mathrm{Frob}_q\circ f$ in which $\mathrm{Frob}_q$ is the Frobenius endomorphism of $X$ and $q$ is a sufficiently large power of $p$.

Now since $q$ is sufficiently large, we see that $g$ is an int-amplified endomorphism of $X$. Moreover, Example 5.1(ii) guarantees that the set $\{n\in\mathbb{N}|\ g^{n}(x)\in Y(K)\}=\{n\in\mathbb{N}|\ f^{n}(x)\in Y(K)\}=\{n\in\mathbb{N}|\ n\cdot P\in V(K)\}$ is \emph{not} a finite union of arithmetic progressions in $\mathbb{N}$.
\end{example}

The Frobenius can even produce polarized dynamical systems which do not satisfy the $\text{DML}_0$ property. We will show this in Example 5.4. It is very interesting that this example comes from the \emph{Drinfeld module}. See \cite[Definition 12.1.2.1]{BGT16}. According to \cite[top of p. 218]{BGT16}, the Denis--Mordell--Lang conjecture about Drinfeld modules motivated Ghioca and Tucker to propose the dynamical Mordell--Lang conjecture.

On the contrary to the convention in \cite[Chapter 12]{BGT16}, the Drinfeld modules in Example 5.4 are of special characteristic. We also recommend the reader to compare Example 5.4 with Corollary 1.9. We notice that variants of Example 5.4 had been studied in \cite{CGLN}.

\begin{example}
Let the base field $K=\overline{\mathbb{F}_{p}(t)}$. Let $\varphi_1(x)=x^p$ and $\varphi_2(y)=y^p+y$. Let $f$ be the split endomorphism of $\mathbb{A}^2$ which sends $(x,y)$ to $(\varphi_1(x),\varphi_2(y))$. In order to get a polarized dynamical system, one can extend $f$ to an endomorphism of either $\mathbb{P}^2$ or $\mathbb{P}^1\times\mathbb{P}^1$. Let $\alpha=(t,t)$ and let the closed subvariety $V$ be the line $x+t=y$. Then $f^n(\alpha)=(t^{p^n},\sum\limits_{i=0}^{n}\binom{n}{i}t^{p^i})$ for every nonnegative integer $n$. Hence $\{n\in\mathbb{N}|\ f^n(\alpha)\in V(K)\}=\{n\in\mathbb{N}|\ t^{p^n}+t=\sum\limits_{i=0}^{n}\binom{n}{i}t^{p^i}\}=\{p^m|\ m\in\mathbb{N}\}$. This gives polarized dynamical systems which do not satisfy the $\text{DML}_0$ property.
\end{example}

\subsection{Endomorphisms with zero entropy}

In this subsection, we will see how complicated the return set can be for endomorphisms with zero entropy. We focus on the endomorphisms of tori because in some sense this is the only case that one can compute the return set. We fix our base field $K=\overline{\mathbb{F}_p(t)}$. In order to be safe, we require our prime $p$ to be not too small, e.g. $p\geq11$.

In $\cite{CGSZ21}$, the authors find that one may reduce the $p$DML problem for tori to the problem of solving polynomial-exponential equations. But we want to emphasize that we need to solve a \emph{system} of polynomial-exponential equations instead of a single one. We start with an heuristic example which corresponds to the system of equations below.
$$
\left\{
\begin{array}{cc}
n=p^{n_1}+p^{n_2} \\
n^2=p^{n_2}+2p^{n_3}+p^{n_4}
\end{array}
\right.
$$

\begin{example}
Consider $f:\mathbb{G}_m^{6}\rightarrow\mathbb{G}_m^{6}$ given by
$$
(x_1,x_2,x_3,x_4,x_5,x_6)\mapsto((t+1)^2x_1,x_1x_2,t^2x_3,x_3x_4,\\(t-1)^2x_5,x_5x_6).
$$
Let the statrting point $\alpha=(t+1,1,t,1,t-1,1)$ and let the closed subvariety $V\subseteq\mathbb{G}_m^{6}$ be $V=\alpha\cdot C_1\cdot C_2\cdot C_3\cdot C_4$ where $C_1,C_2,C_3,C_4\subseteq\mathbb{G}_m^{6}$ are closed subcurves given by
\begin{enumerate}
\item
$C_1=\{((u+1)^2,1,u^2,1,(u-1)^2,1)|\ u\in K\backslash\{0,\pm1\}\}$,
\item
$C_2=\{(v+1)^2,v+1,v^2,v,(v-1)^2,v-1)|\ v\in K\backslash\{0,\pm1\}\}$,
\item
$C_3=\{(1,(w+1)^2,1,w^2,1,(w-1)^2)|\ w\in K\backslash\{0,\pm1\}\}$, and
\item
$C_4=\{(1,x+1,1,x,1,x-1)|\ x\in K\backslash\{0,\pm1\}\}$.
\end{enumerate}
We abuse some notation here since this is just an heuristic example. Now we calculate that $f^n(\alpha)=((t+1)^{2n+1},(t+1)^{n^2},t^{2n+1},t^{n^2},(t-1)^{2n+1},(t-1)^{n^2})$ for every $n\in\mathbb{N}$. Then we get $\{n\in\mathbb{N}|\ f^n(\alpha)\in V(K)\}=\{n\in\mathbb{N}|\ ((t+1)^{2n},(t+1)^{n^2},t^{2n},t^{n^2},(t-1)^{2n},(t-1)^{n^2})\in C_1\cdot C_2\cdot C_3\cdot C_4\}$ and we find that a set of the form $\{p^m+p^{2m}|\ m\in\mathbb{N}\}$ involves in here.
\end{example}

We will use this example to give a rigorous proof of Proposition 1.12 later. Now we give another heuristic example to illustrate that the form of return sets can go beyond the scope of ``widely $p$-normal sets". This example corresponds to the system of equations below.
$$
\left\{
\begin{array}{cc}
n-1=p^{n_1}+p^{n_2}+p^{n_3} \\
n^2-1=2p^{n_1}+2p^{n_2}+4p^{n_3}+2p^{n_4}+2p^{n_5}+p^{n_6}+p^{n_7}+p^{n_8}
\end{array}
\right.
$$

\begin{example}
Consider $f:\mathbb{G}_m^{12}\rightarrow\mathbb{G}_m^{12}$ given by
$$
(x_1,x_2,x_3,x_4,x_5,x_6,x_7,x_8,x_9,x_{10},x_{11},x_{12})\mapsto
$$
$$
(t^2x_1,x_1x_2,(t+1)^2x_3,x_3x_4,(t+2)^2x_5,x_5x_6,(t+3)^2x_7,x_7x_8,(t+4)^2x_9,x_9x_{10},(t+5)^2x_{11},x_{11}x_{12}).
$$
Let the statrting point $\alpha=(t,1,t+1,1,t+2,1,t+3,1,t+4,1,t+5,1)$ and let the closed subvariety $V\subseteq\mathbb{G}_m^{12}$ be $V=\beta\cdot C_1\cdot C_1\cdot C_2\cdot C_3\cdot C_3\cdot C_4\cdot C_4\cdot C_4$ where
$$
\beta=(t^3,t,(t+1)^3,t+1,(t+2)^3,t+2,(t+3)^3,t+3,(t+4)^3,t+4,(t+5)^3,t+5)
$$
and $C_1,C_2,C_3,C_4\subseteq\mathbb{G}_m^{12}$ are closed subcurves given by
\begin{enumerate}
\item
$C_1=\{(u_1^2,u_1^2,(u_1+1)^2,(u_1+1)^2,(u_1+2)^2,(u_1+2)^2,(u_1+3)^2,(u_1+3)^2,(u_1+4)^2,(u_1+4)^2,\\(u_1+5)^2,(u_1+5)^2)|\ u_1\in K\backslash\{0,-1,-2,-3,-4,-5\}\}$,
\item
$C_2=\{(u_2^2,u_2^4,(u_2+1)^2,(u_2+1)^4,(u_2+2)^2,(u_2+2)^4,(u_2+3)^2,(u_2+3)^4,(u_2+4)^2,(u_2+4)^4,\\(u_2+5)^2,(u_2+5)^4)|\ u_2\in K\backslash\{0,-1,-2,-3,-4,-5\}\}$,
\item
$C_3=\{(1,u_3^2,1,(u_3+1)^2,1,(u_3+2)^2,1,(u_3+3)^2,1,(u_3+4)^2,1,(u_3+5)^2)|\ u_3\in K\backslash\{0,-1,-2,\\-3,-4,-5\}\}$, and
\item
$C_4=\{(1,u_4,1,u_4+1,1,u_4+2,1,u_4+3,1,u_4+4,1,u_4+5)|\ u_4\in K\backslash\{0,-1,-2,-3,-4,-5\}\}$.
\end{enumerate}
Calculate as in Example 5.5, we can get $\{n\in\mathbb{N}|\ f^n(\alpha)\in V(K)\}=\{n\in\mathbb{N}|\ (t^{2(n-1)},t^{n^2-1},(t+1)^{2(n-1)},(t+1)^{n^2-1},(t+2)^{2(n-1)},(t+2)^{n^2-1},(t+3)^{2(n-1)},(t+3)^{n^2-1},(t+4)^{2(n-1)},(t+4)^{n^2-1},(t+5)^{2(n-1)},(t+5)^{n^2-1})\in(C_1\cdot C_1\cdot C_2)\cdot(C_3\cdot C_3\cdot C_4\cdot C_4\cdot C_4)\}$ and we find that a set of the form $\{p^{m_1}+p^{m_2}+p^{m_1+m_2}|\ m_1,m_2\in\mathbb{N}\}$ involves in here.
\end{example}

Since Example 5.6 is quite complicated, we shall do a little bit more explanation. The explanation is also heuristic and merely intends to help the reader comprehend this example. The condition on $n$ is a system of 12 equations with 8 variables. We may solve the 6 equations corresponding to the odd coordinates to get the value of the variables $u_1,u_1'$ corresponding to $C_1$ and $u_2$ corresponding to $C_2$. We believe these 6 equations lead to $n-1=p^{n_1}+p^{n_2}+p^{n_3}$ and $(u_1,u_1',u_2)=(t^{p^{n_1}},t^{p^{n_2}},t^{p^{n_3}})$. Then, we may solve another 6 equations corresponding to the even coordinates to get the value of the variables $u_3,u_3'$ corresponding to $C_3$ and $u_4,u_4',u_4''$ corresponding to $C_4$. We believe these 6 equations will somehow lead to $n^2-1=2p^{n_1}+2p^{n_2}+4p^{n_3}+2p^{n_4}+2p^{n_5}+p^{n_6}+p^{n_7}+p^{n_8}$ and $(u_3,u_3',u_4,u_4',u_4'')=(t^{p^{n_4}},t^{p^{n_5}},t^{p^{n_6}},t^{p^{n_7}},t^{p^{n_8}})$. Thus we can calculate the possible values of $n$ by solving those two polynomial-exponential equations. We need the dimension to be that large because we need to guarantee that the number of equations is greater than the number of variables in the procedure above.

One can see that the procedure of translation a system of polynomial-exponential equations into the $p$DML problem of a low-complexity endomorphism of a torus as above is quite free.

~

Now we give a rigorous proof of Proposition 1.12.

\begin{proposition}
Consider the automorphism $f\times g$ of $\mathbb{G}_m^6\times\mathbb{G}_m^3$ in which $f:\mathbb{G}_m^{6}\rightarrow\mathbb{G}_m^{6}$ is defined as in Example 5.5 and $g:\mathbb{G}_m^3\rightarrow\mathbb{G}_m^3$ is the translation $(y_1,y_2,y_3)\mapsto((t+1)y_1,ty_2,(t-1)y_3)$. Let $V_0\subseteq\mathbb{G}_{m}^{3}$ be the closed subvariety given by the equation $y_1+y_3=2y_2+2$ and let $\alpha_0=(1,1,1)$ be the zero element in $\mathbb{G}_{m}^{3}(K)$. Let $\alpha\in\mathbb{G}_{m}^{6}(K)$ be as in Example 5.5. Then there exists a closed subvariety $V\subseteq\mathbb{G}_{m}^{6}$ such that $\{n\in\mathbb{N}|\ (f\times g)^n((\alpha,\alpha_0))\in(V\times V_0)(K)\}$ is \emph{not} a $p$-normal set in $\mathbb{N}$.
\end{proposition}

The proof of Proposition 5.7 is very similar to the proof of \cite[Proposition 5.5]{XY}. We will inherit the notations introduced there. Let $q\in\{p^{n}|\ n\in\mathbb{Z}_+\},q_1,q_2\in\{p^{n}|\ n\in\mathbb{N}\}$ and $c_0\in\mathbb{N},c_1\in\mathbb{Z}_+$. We denote $A(q;q_1,q_2)$ as the set $\{q_1q^{n_1}+q_2q^{n_2}|\ n_1,n_2\in\mathbb{N}\}$ and denote $B(q;c_0,c_1)$ as the set $\{c_0+c_1q^{n}|\ n\in\mathbb{N}\}$. We will obey the convention that all of these coefficients must lie in their ``domain of definition" (i.e. $q\in\{p^{n}|\ n\in\mathbb{Z}_+\},q_1,q_2\in\{p^{n}|\ n\in\mathbb{N}\}$ and $c_0\in\mathbb{N},c_1\in\mathbb{Z}_+$) when we use this notation.

\proof[Proof of Proposition 5.7]
For every closed subvariety $V\subseteq\mathbb{G}_{m}^{6}$, we denote $S(V)$ as the set $\{n\in\mathbb{N}|\ (f\times g)^n((\alpha,\alpha_0))\in(V\times V_0)(K)\}$. Assume by contradiction that $S(V)$ is a $p$-normal set in $\mathbb{N}$ for every closed subvariety $V\subseteq\mathbb{G}_{m}^{6}$. Notice $\{n\in\mathbb{N}|\ g^{n}(\alpha_0)\in V_0(K)\}=\{p^{n_1}+p^{n_2}|\ n_1,n_2\in\mathbb{N}\}$, we may conclude that up to a finite set, $S(V)$ is a union of finitely many sets of the form $A(q;q_1,q_2)$ along with finitely many sets of the form $B(q;c_0,c_1)$ for any closed subvariety $V\subseteq\mathbb{G}_{m}^{6}$ as in the proof of $\cite[\text{Proposition}\ 5.5]{XY}$.

Now let $X$ be the image of the morphism $(\mathbb{A}^1\backslash\{0,\pm1\})^4\rightarrow\mathbb{G}_m^6$ given by
$$
(u,v,w,x)\mapsto
$$
$$
((t+1)(u+1)^2(v+1)^2,(v+1)(w+1)^2(x+1),tu^2v^2,vw^2x,(t-1)(u-1)^2(v-1)^2,(v-1)(w-1)^2(x-1)).
$$
Then $X$ is a constructible set in $\mathbb{G}_m^6$ and hence we may write $X=\bigcup\limits_{i=1}^{N}(V_i\backslash W_i)$ in which $V_1,\dots,V_N,\\W_1,\dots,W_N\subseteq\mathbb{G}_{m}^{6}$ are closed subvarieties satisfying $W_i\subseteq V_i$ for every $1\leq i\leq N$.

So $\bigcup\limits_{i=1}^{N}(S(V_i)\backslash S(W_i))=\{n\in\mathbb{N}|\ g^{n}(\alpha_0)\in V_0(K)\}\cap\{n\in\mathbb{N}|\ f^{n}(\alpha)\in X\}=A(p;1,1)\cap\{n\in\mathbb{N}|\ \exists u,v,w,x\in K\backslash\{0,\pm1\}\ s.t.\ ((t+1)^{2n},(t+1)^{n^2},t^{2n},t^{n^2},(t-1)^{2n},(t-1)^{n^2})=((u+1)^2(v+1)^2,(v+1)(w+1)^2(x+1),u^2v^2,vw^2x,(u-1)^2(v-1)^2,(v-1)(w-1)^2(x-1))\}$. Then we know $\{p^{n}+p^{2n}|\ n\in\mathbb{N}\}\subseteq\bigcup\limits_{i=1}^{N}(S(V_i)\backslash S(W_i))$. Therefore, we can prove that there exist $q_0\in\{p^{n}|\ n\in\mathbb{Z}_+\}$ and $q_{10},q_{20}\in\{p^{n}|\ n\in\mathbb{N}\}$ such that $\{(q_{10}+q_{20})q_0^{n}|\ n\in\mathbb{N}\}\subseteq\bigcup\limits_{i=1}^{N}(S(V_i)\backslash S(W_i))$ as in the proof of $\cite[\text{Proposition}\ 5.5]{XY}$. So we can find $c\in\mathbb{N}$ such that $M=\{n\in\mathbb{N}|\ p^n+p^{n+c}\in\bigcup\limits_{i=1}^{N}(S(V_i)\backslash S(W_i))\}$ is an infinite set.

Now for any $n\in M$, the system of equations
$$
\left\{
\begin{array}{cc}
\pm(t+1)^{p^n+p^{n+c}}=(u+1)(v+1) \\
\pm t^{p^n+p^{n+c}}=uv \\
\pm(t-1)^{p^n+p^{n+c}}=(u-1)(v-1) \\
(t+1)^{p^{2n}(1+p^c)^2}=(v+1)(w+1)^2(x+1) \\
t^{p^{2n}(1+p^c)^2}=vw^2x \\
(t-1)^{p^{2n}(1+p^c)^2}=(v-1)(w-1)^2(x-1)
\end{array}
\right.
$$
has a solution $(u_n,v_n,w_n,x_n)\in(K\backslash\{0,\pm1\})^4$. Since $(u_n+1)(v_n+1)+(u_n-1)(v_n-1)=2u_nv_n+2$, we can see that all of the three ``$\pm$" in the equations must be ``+". So we have $\{u_n,v_n\}=\{t^{p^n},t^{p^{n+c}}\}$ for every $n\in M$. Without loss of generality, we assume that $v_n=t^{p^{n+c}}$ for infinitely many $n$ (the case that $v_n=t^{p^n}$ for infinitely many $n$ can be dealed with by just the same argument as below). Then there is an infinite set $M_1\subseteq M$ such that the system of equations
$$
\left\{
\begin{array}{cc}
(t+1)^{p^{2n}(1+p^c)^2-p^{n+c}}=(w+1)^2(x+1) \\
t^{p^{2n}(1+p^c)^2-p^{n+c}}=w^2x \\
(t-1)^{p^{2n}(1+p^c)^2-p^{n+c}}=(w-1)^2(x-1)
\end{array}
\right.
$$
has a solution $(w_n,x_n)\in(K\backslash\{0,\pm1\})^2$ for every $n\in M_1$. We denote $y_n=t^{p^{n+c}}$ and $z_n=t^{p^{2n}}$ for every $n\in M_1$ and denote $m=(1+p^c)^2$. Then we can see that $\{(y_n,z_n)|\ n\in M_1\}$ is a dense set in $\mathbb{A}_K^2$ and
$$
\left\{
\begin{array}{cc}
\frac{(z_n+1)^m}{y_n+1}=(w_n+1)^2(x_n+1) \\
\frac{z_n^m}{y_n}=w_n^2x_n \\
\frac{(z_n-1)^m}{y_n-1}=(w_n-1)^2(x_n-1)
\end{array}
\right.
$$
So there is an algebraic relation between $\frac{(z_n+1)^m}{y_n+1},\frac{z_n^m}{y_n}$, and $\frac{(z_n-1)^m}{y_n-1}$. Namely, we have
$$
\left(9\cdot\frac{z_n^m}{y_n}-\frac{1}{4}\cdot(\frac{(z_n+1)^m}{y_n+1}-2\cdot\frac{z_n^m}{y_n}+\frac{(z_n-1)^m}{y_n-1})\cdot(\frac{(z_n+1)^m}{y_n+1}-\frac{(z_n-1)^m}{y_n-1}-2)\right)^2=
$$
$$
4\cdot\left(\frac{1}{4}\cdot(\frac{(z_n+1)^m}{y_n+1}-2\cdot\frac{z_n^m}{y_n}+\frac{(z_n-1)^m}{y_n-1})^2-\frac{3}{2}\cdot(\frac{(z_n+1)^m}{y_n+1}-\frac{(z_n-1)^m}{y_n-1}-2)\right)
$$
$$
\cdot\left(\frac{1}{4}\cdot(\frac{(z_n+1)^m}{y_n+1}-\frac{(z_n-1)^m}{y_n-1}-2)^2-\frac{3}{2}\cdot(\frac{(z_n+1)^m}{y_n+1}-2\cdot\frac{z_n^m}{y_n}+\frac{(z_n-1)^m}{y_n-1})\cdot\frac{z_n^m}{y_n}\right)
$$
for every $n\in M_1$. But since $\{(y_n,z_n)|\ n\in M_1\}$ is a dense set in $\mathbb{A}_K^2$, this equation must be an identity with variables $y$ and $z$. However, regarding LHS and RHS as polynomials in $K(y)[z]$, we can calculate that the coefficient of $z^{4m}$ in LHS is $\frac{1}{y^2(y^2-1)^4}$ while this coefficient in RHS is $\frac{4(3-2y^2)}{y^4(y^2-1)^4}$. So this equation cannot be an identity and hence we get a contradiction. Thus we conclude that there exists a closed subvariety $V\subseteq\mathbb{G}_{m}^{6}$ such that $\{n\in\mathbb{N}|\ (f\times g)^n((\alpha,\alpha_0))\in(V\times V_0)(K)\}$ is not a $p$-normal set in $\mathbb{N}$.
\endproof

Proposition 1.12 is an immediate consequence of Proposition 5.7, as all the dynamical systems considered in this subsection have zero entropy.

~

At the end of this article, we formulate a refined version of the $p$DML conjecture, following a referee's suggestion. It seems that except for a little gap in the proof of \cite[Theorem 3.2]{CGSZ21} (which incorrectly believes that solving a system of equations is equivalent to solving each one of them), all of the arguments in that paper are valid. So we tried to follow their methods and find out what can we say about the return sets of endomorphisms of tori after assuming Vojta's conjecture. These thoughts lead to the following question.

We firstly introduce the notion of \emph{semilinear sets} in $\mathbb{N}^n$ \cite{GS66}.

\begin{definition}
Let $n$ be a positive integer. A \emph{linear subset of} $\mathbb{N}^n$ is a coset of a finitely generated sub-semigroup of $\mathbb{N}^n$. A \emph{semilinear subset of} $\mathbb{N}^n$ is a finite union of linear subsets in $\mathbb{N}^n$.
\end{definition}

In $\mathbb{N}^n$, the semilinear subsets are precisely the definable subsets of the Presburger arithmetic of $\mathbb{N}$. Also, this family of sets is closed with respect to intersection, complementation, and projection. See \cite[Section 1]{GS66}.

\begin{question}
Let $X$ be a variety over an algebraically closed field $K$ of characteristic $p>0$. Let $f:X\dashrightarrow X$ be a dominant rational self-map. We introduce the sets
$$
S_{d,\Lambda}(c_{0};c_{1},\dots,c_{d})=\{c_{0}+\sum\limits_{i=1}^{d}c_{i}p^{n_{i}}|\ (n_{1},\dots,n_{d})\in\Lambda\}
$$
in which $d$ is a positive integer, $c_0,c_1,\dots,c_d$ are rational numbers, and $\Lambda\subseteq\mathbb{N}^{d}$ is a semilinear subset. Let $x\in X(K)$ be a point such that the orbit $\mathcal{O}_f(x)$ is well-defined, and let $V\subseteq X$ be a closed subvariety. Then up to a finite set, does the return set $\{n\in\mathbb{N}|\ f^{n}(x)\in V(K)\}$ have to be a finite union of arithmetic progressions in $\mathbb{N}$ along with finitely many sets of the form $S_{d,\Lambda}(c_{0};c_{1},\dots,c_{d})$?
\end{question}

In order to eliminate possible confusions, we remark that in the statement above, we require all sets $S_{d,\Lambda}(c_{0};c_{1},\dots,c_{d})$ involved in the expression to be contained in $\mathbb{N}$.

As we have explained before, the $p$DML problem has a close connection with solving polynomial-exponential equations. In \cite[Theorem 5.5]{CGSZ21}, the authors calculated the solution set of a single polynomial-exponential equation under the assumption of Vojta's conjecture. By looking through its proof, we find that the meaning of ``$p$-sets" in loc. cit. should be comprehended in the sense as Question 5.9 above. Following the method of \cite{CGSZ21} and assuming Vojta's conjecture, maybe one can give a positive answer of this question for endomorphisms of tori. But a careful computation is needed.

\section*{Acknowledgements}
We would like to thank Jiarui Song, Chengyuan Yang, and Geng-Rui Zhang for some useful discussion. We thank the anonymous referees for many beneficial suggestions.

This work is supported by the National Natural Science Foundation of China Grant No. 12271007.

\bibliographystyle{plain}
\bibliography{reference}

\address{Beijing International Center for Mathematical Research, Peking University, Beijing 100871, China}

\email{xiejunyi@bicmr.pku.edu.cn}

~

\address{Beijing International Center for Mathematical Research, Peking University, Beijing 100871, China}

\email{ys-yx@pku.edu.cn}

\end{spacing}
\end{document}